\documentclass[3p]{elsarticle}
\usepackage{bm}
\usepackage[parfill]{parskip}
\usepackage[utf8]{inputenc}
\usepackage{url}
\usepackage{amsmath,amssymb,amsfonts,amsthm}
\usepackage{color}
\usepackage{xcolor}
\usepackage{bm}
\usepackage{graphicx}
\usepackage{subcaption}
\usepackage{lipsum}
\usepackage{amsfonts}
\usepackage{graphicx}
\usepackage{epstopdf}
\usepackage{array}
\usepackage[hidelinks]{hyperref}
\usepackage{lineno}
\usepackage{bbm}
\usepackage{xfrac}
\usepackage{varwidth}
\usepackage[nice]{nicefrac}
\usepackage{yhmath}

\renewcommand{\d}{\,\mathrm{d}}

\usepackage{amsmath}
\usepackage{algorithm}
\usepackage{algpseudocode}

\ifpdf
  \DeclareGraphicsExtensions{.eps,.pdf,.png,.jpg}
\else
  \DeclareGraphicsExtensions{.eps}
\fi

\DeclareMathOperator*{\argmin}{arg\,min}

\makeatletter
\newcommand{\doublehat}[1]{%
\begingroup%
  \let\macc@kerna\z@%
  \let\macc@kernb\z@%
  \let\macc@nucleus\@empty%
  \widehat{\raisebox{.35ex}{\vphantom{\ensuremath{#1}}}\smash{\widehat{#1}}}%
\endgroup%
}
\makeatother

\makeatletter
\newsavebox\myboxA
\newsavebox\myboxB
\newlength\mylenA

\newcommand*\xoverline[2][0.75]{%
    \sbox{\myboxA}{$\m@th#2$}%
    \setbox\myboxB\null
    \ht\myboxB=\ht\myboxA%
    \dp\myboxB=\dp\myboxA%
    \wd\myboxB=#1\wd\myboxA
    \sbox\myboxB{$\m@th\overline{\copy\myboxB}$}
    \setlength\mylenA{\the\wd\myboxA}
    \addtolength\mylenA{-\the\wd\myboxB}%
    \ifdim\wd\myboxB<\wd\myboxA%
       \rlap{\hskip 0.5\mylenA\usebox\myboxB}{\usebox\myboxA}%
    \else
        \hskip -0.5\mylenA\rlap{\usebox\myboxA}{\hskip 0.5\mylenA\usebox\myboxB}%
    \fi}
\makeatother

\usepackage{tikz}
\usetikzlibrary{chains, positioning, arrows.meta, bending, shapes.arrows}
\usetikzlibrary{positioning,shapes,arrows,calc,shapes.geometric,backgrounds,fit}
\usepackage{pifont}
\usetikzlibrary{positioning,shapes,arrows,calc,shapes.geometric,backgrounds,fit}

\tikzstyle{block} = [rectangle,draw,minimum width=2em,align=center,rounded corners, minimum height=2em,scale=1.0]
\tikzstyle{blockleft} = [rectangle,draw,minimum width=2em,align=left,rounded corners, minimum height=2em,scale=1.0]
\tikzstyle{bigblock} = [rectangle,draw,minimum width=8em,align=center,rounded corners, minimum height=4em,scale=1.0]
\tikzstyle{connect} = [draw,-latex']
\tikzstyle{decision} = [diamond, draw, 
    text width=4.5em, text badly centered, node distance=3cm, inner sep=0pt]
\tikzstyle{line} = [draw, -latex']
\tikzstyle{cloud} = [draw, ellipse,fill=red!20, node distance=3cm,
    minimum height=2em]
\tikzstyle{linenoarrow}=[draw]

\makeatletter
\let\save@mathaccent\mathaccent
\newcommand*\if@single[3]{%
  \setbox0\hbox{${\mathaccent"0362{#1}}^H$}%
  \setbox2\hbox{${\mathaccent"0362{\kern0pt#1}}^H$}%
  \ifdim\ht0=\ht2 #3\else #2\fi
  }
\newcommand*\rel@kern[1]{\kern#1\dimexpr\macc@kerna}
\newcommand*\widebar[1]{\@ifnextchar^{{\wide@bar{#1}{0}}}{\wide@bar{#1}{1}}}
\newcommand*\wide@bar[2]{\if@single{#1}{\wide@bar@{#1}{#2}{1}}{\wide@bar@{#1}{#2}{2}}}
\newcommand*\wide@bar@[3]{%
  \begingroup
  \def\mathaccent##1##2{%
    \let\mathaccent\save@mathaccent
    \if#32 \let\macc@nucleus\first@char \fi
    \setbox\z@\hbox{$\macc@style{\macc@nucleus}_{}$}%
    \setbox\tw@\hbox{$\macc@style{\macc@nucleus}{}_{}$}%
    \dimen@\wd\tw@
    \advance\dimen@-\wd\z@
    \divide\dimen@ 3
    \@tempdima\wd\tw@
    \advance\@tempdima-\scriptspace
    \divide\@tempdima 10
    \advance\dimen@-\@tempdima
    \ifdim\dimen@>\z@ \dimen@0pt\fi
    \rel@kern{0.6}\kern-\dimen@
    \if#31
      \overline{\rel@kern{-0.6}\kern\dimen@\macc@nucleus\rel@kern{0.4}\kern\dimen@}%
      \advance\dimen@0.4\dimexpr\macc@kerna
      \let\final@kern#2%
      \ifdim\dimen@<\z@ \let\final@kern1\fi
      \if\final@kern1 \kern-\dimen@\fi
    \else
      \overline{\rel@kern{-0.6}\kern\dimen@#1}%
    \fi
  }%
  \macc@depth\@ne
  \let\math@bgroup\@empty \let\math@egroup\macc@set@skewchar
  \mathsurround\z@ \frozen@everymath{\mathgroup\macc@group\relax}%
  \macc@set@skewchar\relax
  \let\mathaccentV\macc@nested@a
  \if#31
    \macc@nested@a\relax111{#1}%
  \else
    \def\gobble@till@marker##1\endmarker{}%
    \futurelet\first@char\gobble@till@marker#1\endmarker
    \ifcat\noexpand\first@char A\else
      \def\first@char{}%
    \fi
    \macc@nested@a\relax111{\first@char}%
  \fi
  \endgroup
}
\makeatother

\journal{arXiv.org}

\newcommand{\TheTitle}{An adaptive dynamical low-rank optimizer for solving kinetic parameter identification problems}

\date{\today}

\usepackage{amsopn}
\DeclareMathOperator{\diag}{diag}

\journal{arXiv}

\begin{document}
\begin{frontmatter}

\title{\TheTitle}

\author[adressWuerzburg]{Lena Baumann}
\author[adressInnsbruck]{Lukas Einkemmer}
\author[adressWuerzburg]{Christian Klingenberg}
\author[adressAs]{Jonas Kusch}

\address[adressWuerzburg]{University of Wuerzburg, Department of Mathematics,  Wuerzburg, Germany, \href{mailto:lena.baumann@uni-wuerzburg.de}{lena.baumann@uni-wuerzburg.de} (Lena Baumann), \href{mailto:klingen@mathematik.uni-wuerzburg.de}{christian.klingenberg@uni-wuerzburg.de} (Christian Klingenberg) }
\address[adressInnsbruck]{University of Innsbruck, Numerical Analysis and Scientific Computing, Innsbruck, Austria, \href{mailto:lukas.einkemmer@uibk.ac.at}{lukas.einkemmer@uibk.ac.at}}
\address[adressAs]{Norwegian University of Life Sciences, Department of Data Science, \r{A}s, Norway, \href{mailto:jonas.kusch@nmbu.no}{jonas.kusch@nmbu.no}}

\begin{abstract}
The numerical solution of parameter identification problems for kinetic equations can exhibit high computational and memory costs. In this paper, we propose a dynamical low-rank scheme for the reconstruction of the scattering parameters in the radiative transfer equation from a number of macroscopic time-independent measurements. We first introduce a PDE constrained optimization procedure in the continuous setting and derive the adjoint equations using a Lagrangian reformulation. For the scattering coefficient, a periodic B-spline parametrization is used and a gradient descent step for updating its coefficients is formulated. After the discretization, a dynamical low-rank approximation (DLRA) is applied. We make use of the rank-adaptive basis update \& Galerkin integrator to solve both the forward and adjoint problem. We adaptively choose the tolerance of the low-rank solver and the gradient descent step size within a line search approach. Numerical results in two and four dimensions validate the accuracy and efficiency of the proposed DLRA scheme. 
\end{abstract}


\end{frontmatter}

\section{Introduction}\label{sec1:Introduction}

Many practical applications involve non-observable quantities that shall be inferred from related observations and measurements. In medical imaging, for example, a classical problem consists in the non-intrusive reconstruction of properties of an examined tissue from measurements. The propagation of near-infrared light through tissue can be modeled using the \textit{radiative transfer equation (RTE)} \cite{RenBalHielscher2007,KloseNetzBeuthanHielscher2002,KloseHielscher2002}. Neglecting boundary effects, the time-dependent form of this kinetic partial differential equation (PDE) can be given in 1+1D slab geometry as
\begin{align}
\begin{cases}
    \partial_t f\left(t,x,v\right) + v \partial_x f\left(t,x,v\right) &= \sigma \left(x \right) \left(\frac{1}{|\Omega_v|}\langle f\left(t,x,v\right)) \rangle_v - f\left(t,x,v\right) \right),\\
    f\left(t=0,x,v\right) &= f_{\text{in}}\left(x,v\right),
\end{cases}\label{1:RTE}
\end{align}
where $f\left(t,x,v\right) : \mathbb{R}^+_0 \times \Omega_x \times \Omega_v \to \mathbb{R}^+_0$ denotes the distribution function that describes the repartition of photons in phase space. Here, $t$ stands for the time variable, $x \in \Omega_x \subseteq \mathbb{R}$ for the space variable and $v \in \Omega_v= \left[-1,1\right]$ for the angular variable. An integration over the corresponding domain is denoted by brackets $\left\langle \cdot \right\rangle$ and $\left| \Omega_v\right|$ measures the length of the domain $\Omega_v$. The function $\sigma \left(x\right)$ represents the properties of the background medium, indicating the probability of particles at position $x$ to be scattered into a new direction. We refer to it as the \textit{scattering coefficient}. At the initial time $t=0$ the function $f_{\text{in}}\left(x,v\right)$ shall be prescribed for the distribution function. For notational simplicity we focus on the 1+1d case. However, the extension to higher dimensions is immediate and we will present 2+2d (i.e.~four-dimensional) simulations in section \ref{sec:4d-example}.

The parameter identification problem associated to the RTE \eqref{1:RTE} considers the reconstruction of the scattering coefficient $\sigma \left(x\right)$ from measurements. We refer to it as the \textit{inverse transport problem}, for which more general information can be found in the review articles \cite{Bal2009,Stefanov2003}. For the numerical solution of this parameter identification problem, PDE constrained optimization is employed. PDE constrained optimization in the context of kinetic equations has seen some interest lately (see, e.g.~\cite{bartsch2023,einkemmer2024beamheating,albi2025robust,einkemmer2026mirror,lu2025dynamical}). However, this works focus more on optimization and control problems. In our setting, we aim for the minimization of the difference between the measurements and the computed solutions under the assumption of the validity of the RTE. Similar to recent papers \cite{LiWangYang2023,ChenLiLiu2018,HellmuthKlingenbergLiTang2025,EinkemmerLiWangYang2024}, we pursue a gradient-based approach for which in each iteration the evaluation of both the forward and the adjoint problem is required. Clearly, this can numerically become very costly, especially in higher-dimensional settings. 

To reduce the computational cost and memory requirements for the solution of kinetic equations, \textit{dynamical low-rank approximation (DLRA)} \cite{KochLubich2007} can be applied. This approach approximates the kinetic distribution function $f$ up to a certain \textit{rank} $r$ as
\begin{align}\label{1:eq-DLRA}
f \left(t,x,v\right) \approx \sum_{i,j=1}^r X_i \left(t,x \right) S_{ij} \left(t\right) V_j \left(t,v\right),
\end{align}
where $\{X_i : i=1,..,r \}$ are the orthonormal basis functions in space and $\{V_j : j=1,..,r \}$ are the orthonormal basis functions in angle. The matrix $\mathbf{S} = \left(S_{ij}\right) \in \mathbb{R}^{r \times r}$ contains the coefficients of the approximation and therefore is called the \textit{coefficient} or \textit{coupling matrix}. The idea of DLRA then consists in constraining the evolution dynamics to functions of the form \eqref{1:eq-DLRA}. There are integrators that can be used to evolve the low-rank factors in time, while being robust to the presence of small singular values. For instance, the \textit{projector-splitting} \cite{LubichOseledets2014}, the \textit{(augmented) basis update \& Galerkin} (BUG) \cite{CerutiLubich2022,CerutiKuschLubich2022}, and the \textit{parallel BUG integrator} \cite{CerutiKuschLubich2024} are widely used in various areas of research \cite{BaumannEinkemmerKlingenbergKusch2024,KuschStammer2023,EinkemmerLubich2018,EinkemmerHuYing2021}. We also mention the recent review article \cite{lrkin2025review} specifically dedicated to low-rank methods for kinetic equations.

For the solution of the inverse transport problem associated to \eqref{1:RTE}, the following approach is pursued: ``first optimize, then discretize, then low-rank", i.e.~we first perform the optimization in a continuous setting before the resulting equations are discretized and a dynamical low-rank approximation is used. The main features of this paper are: 
\begin{itemize}

    \item \textit{An application of DLRA to a PDE parameter identification problem:} The scattering parameter $\sigma \left( x \right)$ is determined by PDE constrained optimization for which after the discretization the dynamical low-rank method is used. To our knowledge, this is the first paper that combines parameter identification problems and DLRA, leading to a reduction of memory from $\mathcal{O}\left(N^{d_x + d_v}\right)$ to $\mathcal{O}\left(r N^{\max\left(d_x, d_v \right)}\right)$ in each step, where $N$ denotes the number of grid points in physical as well as angular space and $d_x,d_v$ the dimensions in space and angle, respectively. Such savings are particularly important in this context as in order to compute the gradient many problems need to store the entire history of the forward solve.

\item \textit{A setup close to realistic applications:} In most applications measurements are not able to access the full distribution function but at most angle-averaged quantities, i.e. its moments. We will consider such a setup here, where it is assumed that only the first moment is accessible by measurements. In addition, optimal tomography commonly relies on a multitude of measurements from different positions which we incorporate by probing multiple initial values.

\item \textit{An adaptive gradient descent step size and an augmented low-rank integrator:} The minimization is performed using a gradient descent method for updating the coefficients of a periodic B-spline approximation of $\sigma\left(x\right)$. Similar to \cite{ScaloneEinkemmerKuschMcClarren2024}, the step size is adaptively chosen by a backtracking line search approach with Armijo condition. Also the rank of the DLRA algorithm is adaptively determined by using the augmented BUG integrator presented in \cite{CerutiKuschLubich2022} combined with an adaptively chosen truncation error tolerance. As a result, this enables us to begin the optimization procedure with a comparatively small rank (when the solution is still far from the minimum) and gradually increase the rank as the optimization progresses and more accuracy is required, again enhancing the performance of the DLRA scheme.

\item \textit{A series of numerical test examples:} A series of numerical text examples in 1D and 2D confirms that for the reconstruction of the scattering coefficient for the inverse transport problem the application of DLRA shows good agreement with the full solution while being significantly faster, suggesting that the combination of low-rank methods and parameter identification problems is a promising field of future research.

\end{itemize}

The structure of the paper is as follows: After the introduction in Section \ref{sec1:Introduction}, the PDE constrained optimization procedure for the solution of the inverse transport problem is explained in Section \ref{sec2:PDE constrained optimization}. Section \ref{sec3:Discretization} is devoted to the discretization of the forward and the adjoint equations as well as of the gradient in angle, space, and time, leading to a fully discrete gradient descent scheme. In Section \ref{sec4:DLRA} the concept of DLRA is introduced and subsequently applied to the forward and adjoint equations. An adaptive backtracking line search method for an adaptive refinement of the gradient descent step size and the DLRA rank tolerance is presented. Numerical results given in Section \ref{sec5:Numerical results} confirm the accuracy and efficiency of the DLRA scheme compared to the solutions computed with the full solver. Finally, Section \ref{sec6:Conclusion and outlook} gives a brief conclusion and an outlook for possible further research. 

\section{PDE constrained optimization}\label{sec2:PDE constrained optimization}

In practical applications, optical tomography commonly relies on a multitude of measurements from different positions. To be close to realistic settings, we take a number of $N_{\text{IC}}$ measurements for the reconstruction of the scattering coefficient $\sigma\left(x\right)$ into account. Furthermore, we assume the measurements to be generated by a measurement operator $M$ acting on the angle-averaged solution of the RTE at the final time $t=T$, which has been computed with the corresponding initial condition $f_{\text{in},m}$ for $m=1,...,N_{\text{IC}}$. The restriction to final-time measurements is a common choice in PDE constrained optimization \cite{LiWangYang2023,CaflischSilantyevYang2021}, while the assumption of angle-averaged measurements is often made when modeling real-world problems \cite{BalJollivet2009,ChenLiWang2018}. For simplicity, the experimentally observed data $d_m$ is assumed to be close to the measurements of the angle-averaged solution that is obtained
with the ground truth parameter, i.e.
\begin{align*}
d_m \left(x \right) \approx M \left( \left\langle f_{\sigma,m} \left(t=T,x,v \right) \right\rangle_v \right) \quad \text{ for } \quad m=1,...,N_{\text{IC}}, 
\end{align*}
where $f_{\sigma,m} \left(t,x,v\right)$ denotes a solution of
\begin{align}\label{2:RTE multiple IC}
\begin{cases}
\partial_t f_m \left(t,x,v\right) + v \partial_x f_m \left(t,x,v\right) &= \sigma \left(x \right) \left(\frac{1}{|\Omega_v|}\langle f_m \left(t,x,v\right)) \rangle_v - f_m \left(t,x,v\right) \right),\\
f_m\left(t=0,x,v\right) &= f_{\text{in},m}\left(x,v\right),
\end{cases}
\end{align}
computed with the corresponding scattering coefficient $\sigma \left(x \right)$. For the solution of the PDE parameter identification problem associated with \eqref{2:RTE multiple IC} one tries to minimize the square loss between the measurements of the angle-averaged distribution function satisfying \eqref{2:RTE multiple IC} and the experimentally observed data, i.e. one tries to solve the minimization problem
\begin{align}\label{2:Optimization min functional}
\argmin_{\sigma} J\left(\sigma\right) \quad \text{ with } \quad J\left(\sigma \right)= \frac{1}{2} \sum_{m=1}^{N_{\text{IC}}} \left\langle \left| \left\langle M \left(f_{\sigma,m} (t=T,x,v) \right\rangle_v\right) - d_m\left(x\right) \right|^2 \right\rangle_x \quad \text{ subject to } \quad \eqref{2:RTE multiple IC}.
\end{align}

Note that this setup is close to realistic applications in the sense as described above. However, it is assumed here that there is no noise in the measurements which in practical applications is clearly infeasible. Nevertheless, the results gained from the considered setup can directly be extended to higher-dimensional settings and give valuable insights into the combination of parameter identification and DLRA, which is the goal of this paper.

In Section \ref{sec2.1:Lagrangian formulation} we make use of the method of Lagrange multipliers to derive the adjoint equations associated to the forward problem \eqref{2:RTE multiple IC}. We derive the explicit gradient descent step in Section \ref{sec2.2:Gradient descent step}. 

\subsection{Lagrangian formulation and adjoint problem}\label{sec2.1:Lagrangian formulation}

For the derivation of the adjoint problem we reformulate the PDE constrained minimization problem \eqref{2:Optimization min functional} using the method of Lagrange multipliers. We aim for a solution of 
\begin{align*}
\argmin \mathcal{L} \left(f_1,...f_{N_{\text{IC}}}, g_1,...,g_{N_{\text{IC}}},\lambda_1,...,\lambda_{N_{\text{IC}}},\sigma\right),
\end{align*}
where 
\begin{align*}
   \mathcal{L} = J \left(f_1,...,f_{N_{\text{IC}}}, \sigma\right)&+ \sum_{m=1}^{N_{\text{IC}}} \left\langle g_m, \partial_t f_m + v \partial_x f_m - \sigma \left(x \right) \left( \frac{1}{|\Omega_v|} \left\langle f_m \right\rangle_v -f_m \right)\right\rangle_{t,x,v}\\
   &+ \sum_{m=1}^{N_{\text{IC}}}\left\langle \lambda_m, f_m \left(t=0,x,v \right) - f_{\text{in},m} (x,v) \right\rangle_{x,v}
\end{align*}
and $g_m \left(t,x,v\right)$ and $\lambda_m \left(x,v\right)$ are the \textit{Lagrange multipliers} with respect to $f_m\left(t,x,v\right)$ and $f_{\text{in},m} \left(x,v\right)$ for $m=1,...,N_{\text{IC}}$, respectively. Applying integration by parts and assuming periodic boundary conditions, the Lagrangian can be rewritten as
\begin{align*}
\mathcal{L} = J \left(f_1,...,f_{N_{\text{IC}}},\sigma\right) &+ \sum_{m=1}^{N_{\text{IC}}} \left\langle f_m, - \partial_t g_m - v \partial_x g_m - \sigma \left(x\right) \left(\frac{1}{|\Omega_v|} \left\langle g_m \right\rangle_v -g_m \right)\right\rangle_{t,x,v}\\
&+ \sum_{m=1}^{N_{\text{IC}}} \left\langle g_m(t=T,x,v), f_m(t=T,x,v)\right\rangle_{x,v} - \sum_{m=1}^{N_{\text{IC}}} \left\langle g_m \left(t=0,x,v \right), f_m\left( t=0,x,v \right)\right\rangle_{x,v}\\
&+ \sum_{m=1}^{N_{\text{IC}}} \left\langle \lambda_m, f_m \left(t=0,x,v \right) - f_{\text{in},m} \left(x,v \right) \right\rangle_{x,v}.
\end{align*}
The corresponding \textit{adjoint} or \textit{dual problems} associated with \eqref{2:RTE multiple IC} can be derived by setting $\frac{\partial \mathcal{L}}{\partial f_m} = 0$ for $m=1,...,N_{\text{IC}}$. By straightforward calculation one obtains
\begin{align}\label{2:adjoint equations}
\begin{cases}
-\partial_t g_m \left(t,x,v \right) - v \partial_x g_m\left(t,x,v\right) &= \sigma \left(x\right) \left(\frac{1}{|\Omega_v|} \left\langle g_m \left(t,x,v\right) \right\rangle_v -g_m\left(t,x,v\right) \right),\\
g_m \left(t=T,x,v \right) &= - \left\langle f_m\left(t=T,x,v\right) \right\rangle_v+ d_m\left(x\right).
\end{cases}
\end{align}
The notation $g_{\sigma,m} \left(t,x,v \right)$ indicates that $g_m \left(t,x,v \right)$ fulfills equations \eqref{2:adjoint equations}. On the solution manifold with $f_{\sigma,m} \left(t,x,v \right)$ satisfying the PDE constraints \eqref{2:RTE multiple IC} and $g_{\sigma,m} \left(t,x,v \right)$ satisfying the adjoint equations \eqref{2:adjoint equations} for $m=1,...,N_{\text{IC}}$, the following equality holds
\begin{align*}
\mathcal{L}\left(f_{\sigma,1},..., f_{\sigma,{N_{\text{IC}}}}, g_{\sigma,1},...,g_{\sigma,{N_{\text{IC}}}} ,\lambda_1,...,\lambda_{N_{\text{IC}}}, \sigma \right) = J \left( f_{\sigma,1},...,f_{\sigma,{N_{\text{IC}}}} ,\sigma\right),
\end{align*}
translating to the derivatives such that
\begin{align}\label{2:gradientwithL}
\frac{\d J \left( f_{\sigma,1},..., f_{\sigma,{N_{\text{IC}}}},\sigma \right)}{\d \sigma} =& \ \frac{\d\mathcal{L} \left(f_{\sigma,1},..., f_{\sigma,{N_{\text{IC}}}}, g_{\sigma,1},...,g_{\sigma,{N_{\text{IC}}}} ,\lambda_1,...,\lambda_{N_{\text{IC}}}, \sigma \right)}{\d \sigma}\\
=& \ \sum_{m=1}^{N_{\text{IC}}} \left(\frac{\partial \mathcal{L}}{\partial f_{\sigma,m}} \frac{\partial f_{\sigma,m}}{\partial \sigma} + \frac{\partial \mathcal{L}}{\partial g_{\sigma,m}} \frac{\partial g_{\sigma,m}}{\partial \sigma} + \frac{\partial \mathcal{L}}{\partial \lambda_m} \frac{\partial \lambda_m}{\partial \sigma} + \frac{\partial \mathcal{L}}{\partial \sigma} \right) = \sum_{m=1}^{N_{\text{IC}}} \frac{\partial \mathcal{L}}{\partial \sigma}, 
\end{align}
where the first three terms vanish since \eqref{2:RTE multiple IC} and \eqref{2:adjoint equations} are fulfilled. This can be used for an efficient computation of the gradient in the following gradient descent step.

\subsection{Optimization parameters and gradient descent step}\label{sec2.2:Gradient descent step}

\begin{figure}
\centering
\begin{minipage}{0.4\textwidth}
\centering
\includegraphics[width=0.9\linewidth]{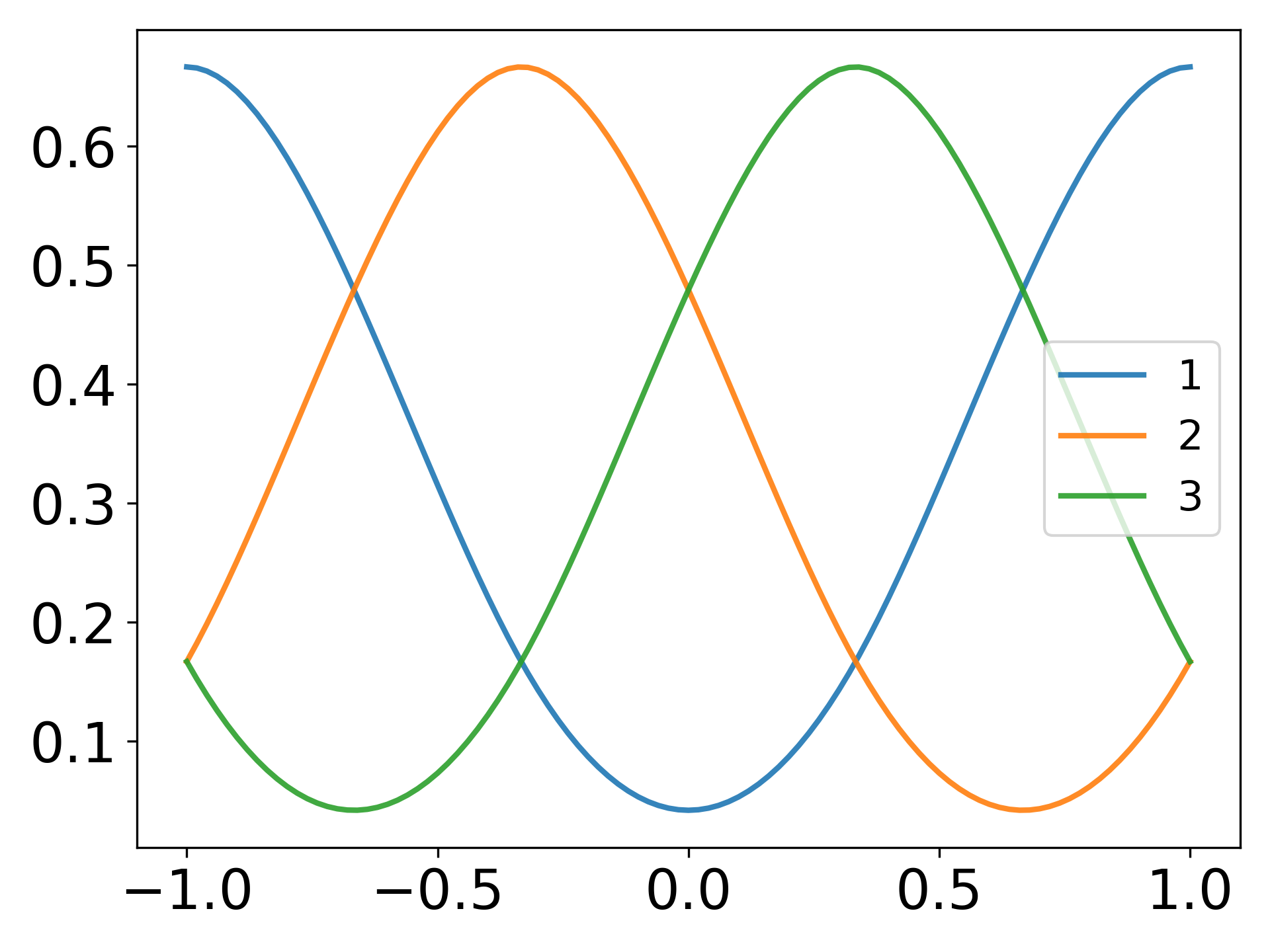}
\end{minipage}
\begin{minipage}{0.4\textwidth}
\centering
\includegraphics[width=0.9\linewidth]{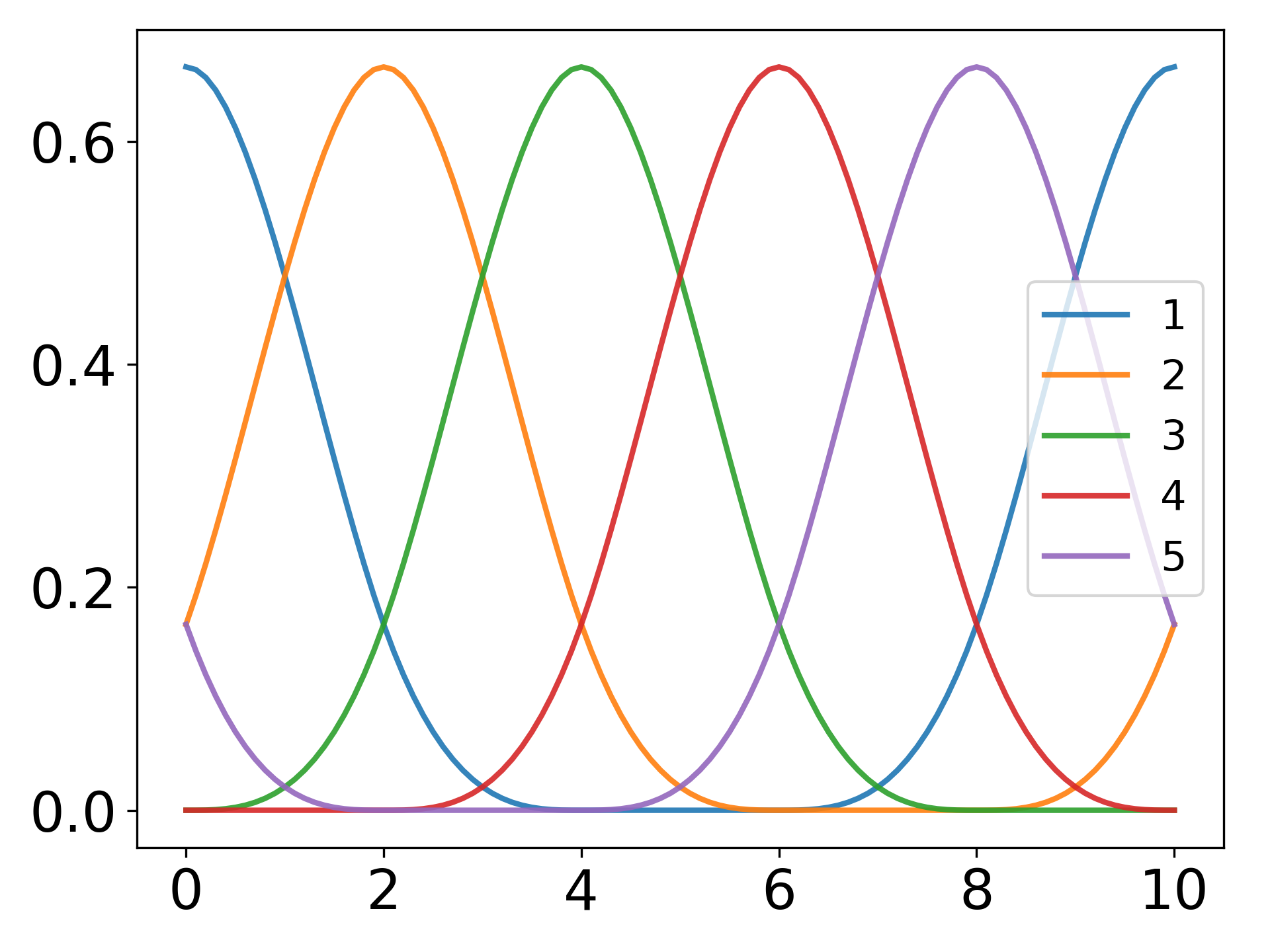}
\end{minipage}
\caption{Cubic periodic B-spline basis functions for $N_c=3$ (left) and $N_c=5$ (right) on different spatial domains.}
\label{fig:Spline basis}
\end{figure}

To avoid computational disadvantages from a large optimization parameter space, the function $\sigma \left(x\right)$ is approximated using splines. In particular, we approximate
\begin{align}\label{2:sigma spline approximation}
\sigma \left(x\right) \approx \sum_{i=1}^{N_c} c_i B_i\left(x\right) \qquad \text{ with } \mathbf{c} = \left(c_i \right) \in \mathbb{R}^{N_c},
\end{align}
where $N_c$ denotes a finite number of cubic periodic B-spline basis functions $B_i \left( x\right)$ and the vector $\mathbf{c}$ contains the coefficients of the approximation. In Figure \ref{fig:Spline basis} the basis functions for the cubic periodic B-splines for $N_c=3$ and $N_c=5$ are illustrated. 

With prescribed initial values $\mathbf{c}^0 \in \mathbb{R}^{N_c}$ the gradient descent step for the solution of the minimization problem \eqref{2:Optimization min functional} updates the coefficients $\mathbf{c}^n$ to $\mathbf{c}^{n+1}$ in each step by determining
\begin{align}\label{2:gradient descent step}
\mathbf{c}^{n+1} = \mathbf{c}^n - \eta^n \nabla_{\mathbf{c}} J \left( \mathbf{c}^n\right) \qquad \text{ for } n=0,1,...,
\end{align}
where $\eta^n$ denotes an adaptively chosen step size. From expression \eqref{2:gradientwithL} we can derive an explicit formulation for the components of the gradient of the cost function as
\begin{align}\label{2:gradient}
\frac{\partial J \left( \mathbf{c} \right)}{\partial c_i} = \sum_{m=1}^{N_{\text{IC}}} \frac{\partial \mathcal{L}}{\partial \sigma} \frac{\partial \sigma}{\partial c_i} = \sum_{m=1}^{N_{\text{IC}}} \left(-\frac{1}{|\Omega_v|} \left\langle \langle f_{\sigma,m} \rangle_v, \langle g_{\sigma,m} \rangle_v \right \rangle_t + \left \langle f_{\sigma,m}, g_{\sigma,m} \right\rangle_{t,v}\right) B_i, 
\end{align}
depending on the solutions of the forward and the adjoint equation as well as on the B-spline basis functions and allowing for an efficient computation of the gradient in the gradient descent step \eqref{2:gradient descent step}. Note that from now on we write $f_m \left(t,x,v\right)$ instead of $f_{\sigma,m} \left(t,x,v\right)$ and $g_m \left(t,x,v\right)$ instead of $g_{\sigma,m} \left(t,x,v\right)$.

\section{Discretization}\label{sec3:Discretization}

For the numerical implementation we discretize the forward problem \eqref{2:RTE multiple IC}, the adjoint problem \eqref{2:adjoint equations} and the gradient \eqref{2:gradient} in angle, space and time, leading to a fully discrete scheme. We begin with the angular discretization in Section \ref{sec3.1:Angular discretization}, followed by a discretization in space in Section \ref{sec3.2:Spatial discretization} and in time in Section \ref{sec3.3:Time discretization}. Section \ref{sec3.4:Fully discrete scheme} summarizes the fully discrete gradient descent method. 

\subsection{Angular discretization}\label{sec3.1:Angular discretization}

For the discretization in angle we decide on a modal approach making use of normalized Legendre polynomials $P_\ell$. This is a standard approach which is commonly used for radiative transfer problems \cite{Pomraning1973,CaseZweifel1967,McClarrenHollowayBrunner2008-2} and has also already been adapted in the inverse setting \cite{WrightSchweigerArridge2007}. We use a rescaling of the Legendre polynomials such that $\left\langle P_k, P_\ell \right\rangle_v = \delta_{k\ell}$ and $P_0 = \frac{1}{\sqrt{2}}$ holds. They constitute a complete set of orthonormal functions on the interval $\left[-1,1 \right]$. We expand the distribution functions $f_m \left(t,x,v \right)$ and $g_m \left(t,x,v \right)$ for $m=1,...,N_{\text{IC}}$ in terms of the rescaled Legendre polynomials and obtain the following approximations
\begin{align}\label{3:angular representation of f and g}
f_m \left( t,x,v \right) \approx \sum_{\ell= 0}^{N_v-1} u_{\ell m} \left( t,x \right) P_\ell \left(v \right) \quad \text{ and } \quad g_m \left( t,x,v \right) \approx \sum_{\ell= 0}^{N_v-1} w_{\ell m} \left( t,x \right) P_\ell \left(v \right),
\end{align}
where the expansion coefficients $u_\ell \left( t,x \right)$ and $w_\ell \left( t,x \right)$, respectively, are called the \textit{moments} and $N_v$ is called the \textit{order} of the approximation. We insert these representations into the forward problem \eqref{2:RTE multiple IC} as well as the adjoint problem \eqref{2:adjoint equations}, multiply with $P_k \left(v\right)$ and integrate over the angular variable $v$. Using the orthonormality condition from above and the notation $A_{k\ell} = \left \langle P_k, v P_\ell \right \rangle_v$ leads to 
\begin{align}\label{3:angular discretization forward equations}
\begin{cases}
\partial_t u_{km} \left( t,x \right) &= - \sum_{\ell=0}^{N_v-1} \partial_x u_{\ell m} \left( t,x \right) A_{k\ell} + \sigma \left(x \right) u_{km} \left( t,x \right) \left(\delta_{k0} - 1 \right),\\
u_{km} \left(t=0,x \right) &= u_{\text{in},km} \left(x\right),
\end{cases}
\end{align}
for the forward equations and
\begin{align}\label{3:angular discretization adjoint equations}
\begin{cases}
- \partial_t w_{km} \left( t,x \right) &= \sum_{\ell=0}^{N_v-1} \partial_x w_{\ell m} \left( t,x \right) A_{k\ell} + \sigma\left(x\right) w_{km} \left( t,x \right) \left(\delta_{k0} - 1 \right),\\
w_{km} \left( t=T,x \right) &= \left(-2 u_{0m} \left(t=T,x\right) + \sqrt{2} d_m \left(x\right) \right) \delta_{k0}, 
\end{cases}
\end{align}
for the adjoint equations for $m=1,...,N_{\text{IC}}$. We collect the entries $A_{k\ell}$ in the symmetric matrix $\mathbf{A} = \left( A_{k\ell}\right) \in \mathbb{R}^{N_v \times N_v}$ and note that $\mathbf{A}$ is diagonalizable in the form $\mathbf{A} = \mathbf{Q}\mathbf{M}\mathbf{Q}^\top$ with $\mathbf{Q}$ orthonormal, $\mathbf{M} = \diag(\sigma_1,...,\sigma_{N_v})$ and $\mathbf{Q}^\top$ denotes the transpose of $\mathbf{Q}$. Then we set $|\mathbf{A}| = \mathbf{Q} |\mathbf{M}| \mathbf{Q}^\top$, where $|\mathbf{M}| = \diag \left( |\sigma_1|,...,|\sigma_{N_v}| \right)$. For the angular discretization of the gradient we insert the representations \eqref{3:angular representation of f and g} into \eqref{2:gradient} and obtain
\begin{align}\label{3:angular discretization gradient}
\frac{\partial J \left( \mathbf{c} \right)}{\partial c_i} \approx \sum_{m=1}^{N_{\text{IC}}} \left(- \left\langle u_{0m} \left(t,x\right), w_{0m} \left(t,x\right) \right \rangle_t + \sum_{k=0}^{N_v-1} \left \langle u_{km} \left(t,x\right), w_{km} \left(t,x\right) \right\rangle_t \right) B_i \left(x\right).
\end{align}

\subsection{Spatial discretization}\label{sec3.2:Spatial discretization}

The discretization in the spatial variable is performed on a spatial grid with $N_x$ grid cells and equidistant spacing $\Delta x = \frac{1}{N_x}$. Spatially dependent quantities are approximated at
the grid points $x_j$ for $j=1,...,N_x$ and denoted by
\begin{align*}
u_{jkm}\left(t\right) \approx u_{km}\left(t,x_j\right), \quad w_{jkm}\left(t\right) \approx w_{km}\left(t,x_j\right), \quad \sigma_j \approx \sigma \left( x_j \right), \quad d_{jm} \approx d_m \left( x_j \right), \quad B_{ji} \approx B_i \left( x_j \right).   
\end{align*}
First-order spatial derivatives are approximated using a centered finite difference scheme to which a second-order stabilization term is added. This discretization is inspired by the stability considerations performed in \cite{BaumannEinkemmerKlingenbergKusch2024}. We denote $\partial_x \approx \mathbf{D}^x \in \mathbb{R}^{N_x \times N_x}$ and $\partial_{xx} \approx \mathbf{D}^{xx} \in \mathbb{R}^{N_x \times N_x}$ for the tridiagonal stencil matrices with nonzero entries only at 
\begin{align*}
D_{j,j\pm 1}^{x}= \frac{\pm 1}{2\Delta x}\;,\qquad D_{j,j}^{xx}= - \frac{2}{\left(\Delta x\right)^2}\;, \quad D_{j,j\pm 1}^{xx}= \frac{1}{\left(\Delta x\right)^2}\;.
\end{align*}
In addition, we assume periodic boundary conditions which results in setting
\begin{align*}
D_{1,N_x}^{x} = \frac{-1}{2\Delta x}
\;,\qquad D_{N_x,1}^{x} = \frac{1}{2\Delta x}\;, \qquad D_{1,N_x}^{xx} = D_{N_x,1}^{xx} = \frac{1}{\left( \Delta x \right)^2}\;.
\end{align*}
The spatially discretized forward equations with centered finite differences and additional second-order stabilization term are obtained from \eqref{3:angular discretization forward equations} as
\begin{align}\label{3:spatial discretization forward problem}
\begin{cases}
\partial_t u_{jkm} \left( t\right) &= - \sum_{i=1}^{N_x} \sum_{\ell=0}^{N_v-1}  D_{ji}^x u_{i\ell m} \left( t \right) A_{k\ell} + \frac{\Delta x}{2} \sum_{i=1}^{N_x} \sum_{\ell=0}^{N_v-1}  D_{ji}^{xx} u_{i\ell m} \left( t \right) \left| A \right|_{k\ell}\\
&\quad  + \ \sigma_j u_{jkm} \left( t \right) \left(\delta_{k0} - 1 \right),\\
u_{jkm} \left(t=0 \right) &= u_{\text{in},jkm},
\end{cases}
\end{align}
and the spatially discretized adjoint equations from \eqref{3:angular discretization adjoint equations} as
\begin{align}\label{3:spatial discretization adjoint problem}
\begin{cases}
- \partial_t w_{jkm} \left( t\right) &= \sum_{i=1}^{N_x} \sum_{\ell=0}^{N_v-1} D_{ji}^x w_{i \ell m} \left( t \right) A_{k\ell} + \frac{\Delta x}{2} \sum_{i=1}^{N_x} \sum_{\ell=0}^{N_v-1} D_{ji}^{xx} w_{i \ell m} \left( t \right) \left|A \right|_{k\ell}\\
&\quad + \ \sigma_j w_{jkm} \left( t\right) \left(\delta_{k0} - 1 \right),\\
w_{jkm} \left( t=T \right) &= \left(-2 u_{j0m} \left(t=T\right) + \sqrt{2} d_{jm} \right) \delta_{k0}. 
\end{cases}
\end{align}
Using the expression given in \eqref{3:angular discretization gradient}, we derive the equation
\begin{align}\label{3:spatial discretization gradient}
\frac{\partial J \left( f_1,...,f_{N_{\text{IC}}} \right)}{\partial c_i} \approx \sum_{m=1}^{N_{\text{IC}}} \left(- \left\langle u_{j0m} \left(t\right), w_{j0m} \left(t\right) \right \rangle_t + \sum_{k=0}^{N_v-1} \left \langle u_{jkm} \left(t\right), w_{jkm} \left(t\right) \right\rangle_t \right) B_{ji}
\end{align}
for the spatial discretization of the gradient of the cost function. The spatially discretized scattering coefficient $\boldsymbol{\sigma} = \left( \sigma_j \right) \in \mathbb{R}^{N_x}$ can be computed as 
\begin{align}\label{2:sigma spline approximation discretized}
\sigma_j \approx \sum_{i=1}^{N_c} c_i B_{ji}.
\end{align}

\subsection{Time discretization}\label{sec3.3:Time discretization}

To obtain a fully discrete system, the time interval $\left[0,T\right]$ is equidistantly split into a finite number $N_t$ of time cells. An update of the forward equations \eqref{3:spatial discretization forward problem} from time $t_n$ to time $t_{n+1} = t_n + \Delta t$ is computed using an explicit Euler step forward in time according to 
\begin{align}\label{3:time discretization forward problem}
\begin{cases}
u_{jkm}^{n+1} &= u_{jkm}^n - \Delta t \sum_{i=1}^{N_x} \sum_{\ell=0}^{N_v-1}  D_{ji}^x u_{i\ell m}^n A_{k\ell} + \Delta t \frac{\Delta x}{2} \sum_{i=1}^{N_x} \sum_{\ell=0}^{N_v-1}  D_{ji}^{xx} u_{i\ell m}^n \left| A \right|_{k\ell}\\
&\quad  + \ \sigma_j \Delta t u_{jkm}^n \left(\delta_{k0} - 1 \right),\\
u_{jkm}^0 &= u_{\text{in},jkm}.
\end{cases}
\end{align}
For the adjoint equations \eqref{3:spatial discretization adjoint problem} we start computations with an end time condition after $N_t$ steps and evolve the solution from time $t_{n}$ to time $t_{n-1} = t_{n} -\Delta t$ by an explicit Euler step backwards in time according to
\begin{align}\label{3:time discretization adjoint problem}
\begin{cases}
w_{jkm}^{n-1} &= w_{jkm}^n + \Delta t \sum_{i=1}^{N_x} \sum_{\ell=0}^{N_v-1} D_{ji}^x w_{i \ell m}^n A_{k\ell} + \Delta t \frac{\Delta x}{2} \sum_{i=1}^{N_x} \sum_{\ell=0}^{N_v-1} D_{ji}^{xx} w_{i \ell m}^n \left|A \right|_{k\ell}\\
&\quad + \ \sigma_j \Delta t w_{jkm}^n \left(\delta_{k0} - 1 \right),\\
w_{jkm}^{N_t} &= \left(-2 u_{j0m}^{N_t}+ \sqrt{2} d_{jm} \right) \delta_{k0}. 
\end{cases}
\end{align}
The fully discretized gradient of the cost function can be obtained from \eqref{3:spatial discretization gradient} by approximating integrals with respect to time by step functions, which yields
\begin{align}\label{3:time discretization gradient}
\frac{\partial J \left( \mathbf{c}\right)}{\partial c_i} \approx  \frac{1}{N_t+1} \sum_{m=1}^{N_{\text{IC}}} \sum_{n=0}^{N_t} \left(- u_{j0m}^n w_{j0m}^n + \sum_{k=0}^{N_v-1} u_{jkm}^n w_{jkm}^n \right) B_{ji}.
\end{align}

\subsection{Fully discrete optimization scheme}\label{sec3.4:Fully discrete scheme}

The strategy for the fully discrete gradient descent method for the solution of the PDE parameter identification problem is summarized in Algorithm \ref{alg:GD fully discrete}. Note that for the stopping criterion an error estimate \texttt{estimated-err} for the deviation of the computed coefficients from the true coefficients is required to run the algorithm. In practical applications, a stopping criterion depending on the amount of progress that is still made in the optimization procedure can be used.

\begin{algorithm}[h!]\small
\caption{\small Gradient descent method for the PDE parameter identification}\label{alg:GD fully discrete}
\begin{algorithmic}
\Require \hspace{0.18cm} measurements $\mathbf{d}_m = \left( d_{jm} \right) \in \mathbb{R}^{N_x}$ for $m=1,...,N_{\text{IC}}$,\\
\hspace{1.04cm} initial data $\mathbf{u}_m^0 = \left(u_{jkm}^0 \right) \in \mathbb{R}^{N_x \times N_v}$ for $m=1,...,N_{\text{IC}}$,\\
\hspace{1.04cm} initial guess for the coefficients $\mathbf{c}^0 = \left(c_i^0 \right)\in \mathbb{R}^{N_c}$,\\
\hspace{1.04cm} initial step size $\eta^0 > 0$,\\ 
\hspace{1.04cm} estimated error \texttt{estimated-err},\\
\hspace{1.04cm} error tolerance \texttt{errtol},\\
\hspace{1.04cm} maximal number of iterations \texttt{maxiter}

\Ensure optimal coefficients $\mathbf{c}^{\text{opt}} = \left(c^{\text{opt}}_i\right) \in \mathbb{R}^{N_c}$ within the prescribed error tolerance\\

\While{\texttt{estimated-err} $>$ \texttt{errtol} \textbf{and} $n\leq$ \texttt{maxiter}}
\State Compute $\bm{\sigma}^n = \left( \sigma^n_j\right) \in \mathbb{R}^{N_x}$ from the given coefficients $\mathbf{c}^n$ according to \eqref{2:sigma spline approximation};
\State Solve the forward problem according to \eqref{3:time discretization forward problem} for each $m=1,...,N_{\text{IC}}$;
\State Solve the adjoint problem according to \eqref{3:time discretization adjoint problem} for each $m=1,...,N_{\text{IC}}$;
\State Compute the components of the gradient $\frac{\partial J \left( \mathbf{c}^n \right)}{\partial c_i^n}$ using \eqref{3:time discretization gradient} and the solutions of \eqref{3:time discretization forward problem} and \eqref{3:time discretization adjoint problem};
\State Update the coefficients according to \eqref{2:gradient descent step}: $\mathbf{c}^{n+1} = \mathbf{c}^n - \eta^n \nabla_{\mathbf{c}} J \left( \mathbf{c}^n\right)$, where $\eta^n$ is adaptively determined by line search;
\EndWhile
\end{algorithmic}
\end{algorithm}\normalsize

\section{Dynamical low-rank approximation}\label{sec4:DLRA}

For the solution of the PDE parameter identification problem the coefficients $\mathbf{c}$ of the spline approximation \eqref{2:sigma spline approximation} of the scattering coefficient $\sigma$ are updated several times in the gradient descent step \eqref{2:gradient descent step}. For each iteration the solution of the fully discretized forward equations \eqref{3:time discretization forward problem} as well as of the fully discretized adjoint equations \eqref{3:time discretization adjoint problem} have to be computed and stored in order to compute the fully discretized gradient \eqref{3:time discretization gradient}. A method for the reduction of computational and memory effort for kinetic equations is the concept of dynamical low-rank approximation, which shall be applied to the considered inverse transport problem. We begin with some general information on DLRA in Section \ref{sec4.1:Background on DLRA}, before Section \ref{sec4.2:DLRA for the optimization} is devoted to a DLRA algorithm for the considered discrete optimization problem.

\subsection{Background on dynamical low-rank approximation}\label{sec4.1:Background on DLRA}

In \cite{KochLubich2007} the concept of DLRA has been introduced in a semi-discrete time-dependent matrix setting. We follow the explanations there. Let $\mathbf{f}\left(t\right) \in \mathbb{R}^{N_x \times N_v}$ be the solution of the matrix differential equation
\begin{align*}
\dot{\mathbf{f}}\left(t\right) = \mathbf{F} \left(\mathbf{f}\left(t\right) \right),
\end{align*}
for which the right-hand side shall be denoted by $\mathbf{F} \left(\mathbf{f}\left(t\right) \right): \mathbb{R}^{N_x \times N_v} \to \mathbb{R}^{N_x \times N_v}$. We seek an approximation of $\mathbf{f}\left(t\right)$ of the form
\begin{align}\label{4:DLRA f=XSV} \mathbf{f}_r \left(t\right) = \mathbf{X}\left(t\right) \mathbf{S}\left(t\right) \mathbf{V}\left(t\right)^\top,
\end{align}
where the matrix $\mathbf{X}\left(t\right) \in \mathbb{R}^{N_x \times r}$ contains the orthonormal basis functions in space and $\mathbf{V}\left(t\right) \in \mathbb{R}^{N_v \times r}$ the orthonormal basis functions in angle. The coefficients of the approximation are stored in the coupling matrix $\mathbf{S}\left(t\right) \in \mathbb{R}^{r \times r}$. The set of all matrices of the form \eqref{4:DLRA f=XSV} constitutes a low-rank manifold that we denote by $\mathcal{M}_r$. Its tangent space at $\mathbf{f}_r\left(t\right)$ shall be denoted by $\mathcal{T}_{\mathbf{f}_r\left(t\right)} \mathcal{M}_r$. For the evolution of the low-rank factors in time we seek a solution of the minimization problem
\begin{align*}
\min_{\dot{\mathbf{f}}_r\left(t\right) \in \mathcal{T}_{\mathbf{f}_r\left(t\right)} \mathcal{M}_r} \left \Vert \dot{\mathbf{f}}_r \left(t\right) - \mathbf{F}\left(\mathbf{f}_r \left(t\right)\right) \right\Vert_F
\end{align*}
at all times $t$, where $\left\Vert \cdot \right\Vert_F $ denotes the Frobenius norm. In \cite{KochLubich2007} it has been shown that this minimization constraint is equivalent to determining 
\begin{align}\label{4:DLRA projection}
\dot{\mathbf{f}}_r\left(t\right) = \mathbf{P} \left(\mathbf{f}_r\left(t\right) \right) \mathbf{F}\left(\mathbf{f}_r \left(t\right)\right), 
\end{align}
where $\mathbf{P}$ denotes the orthogonal projector onto the tangent space $\mathcal{T}_{\mathbf{f}_r\left(t\right)} \mathcal{M}_r$, which can be explicitly given as 
\begin{align*}
\mathbf{P}\left(\mathbf{f}_r \left(t\right) \right) \mathbf{F} = \mathbf{X} \mathbf{X}^\top \mathbf{F} - \mathbf{X} \mathbf{X}^\top \mathbf{F} \mathbf{V} \mathbf{V}^\top + \mathbf{F} \mathbf{V} \mathbf{V}^\top.
\end{align*}
There are different robust time integrators for the solution of \eqref{4:DLRA projection} that are able to evolve the low-rank solution on the manifold $\mathcal{M}_r$ while not suffering from potentially small singular values \cite{KieriLubichWalach2016}. The projector-splitting \cite{LubichOseledets2014}, the (augmented) BUG \cite{CerutiLubich2022,CerutiKuschLubich2022}, and the parallel BUG integrator \cite{CerutiKuschLubich2024} are frequently used.

In this work, we use the augmented BUG integrator from \cite{CerutiKuschLubich2022} which evolves the low-rank factors as follows: In the first two steps, the BUG integrator updates and augments the spatial basis $\mathbf{X}$ and the angular basis $\mathbf{V}$ in parallel, leading to an increase of rank from $r$ to $2r$. Having the augmented bases at hand, a Galerkin step for the coefficient matrix $\mathbf{S}$ is performed. In the last step, all quantities are truncated back to a new rank $r_1 \leq 2r$ that is chosen adaptively depending on a prescribed error tolerance. In detail, the augmented BUG integrator evolves the low-rank solution from $\mathbf{f}_r^n = \mathbf{X}^n \mathbf{S}^n \mathbf{V}^{n,\top}$ at time $t_n$ to $\mathbf{f}_r^{n+1} = \mathbf{X}^{n+1} \mathbf{S}^{n+1} \mathbf{V}^{n+1,\top}$ at time $t_{n+1} = t_n + \Delta t$ as follows: 

\textbf{\textit{K}-Step}: We denote $\mathbf{K}\left(t\right) = \mathbf{X}\left(t\right) \mathbf{S}\left(t\right)$ and solve the PDE 
\begin{align*}
\dot{\mathbf{K}}\left(t\right) = \mathbf{F}\left( \mathbf{K}\left(t\right) \mathbf{V}^{n,\top} \right) \mathbf{V}^n, \quad \mathbf{K}\left(t_n\right) = \mathbf{X}^n \mathbf{S}^n.
\end{align*}
The spatial basis is updated by determining $\widehat{\mathbf{X}}^{n+1} \in \mathbb{R}^{N_x \times 2r}$ as an orthonormal basis of $[\mathbf{K}(t_{n+1}), \mathbf{X}^n] \in \mathbb{R}^{N_x \times 2r}$, e.g. by QR-decomposition. We store $\widehat{\mathbf{M}} = \widehat{\mathbf{X}}^{n+1,\top} \mathbf{X}^n \in \mathbb{R}^{2r \times r}$. Note that we denote augmented quantities of rank $2r$ with hats.

\textbf{\textit{L}-Step}:
We denote $\mathbf{L}\left(t\right) = \mathbf{V}\left(t\right) \mathbf{S}(t)^\top$ and solve the PDE
\begin{align*}
\dot{\mathbf{L}}\left(t\right) = \mathbf{F}\left(\mathbf{X}^n \mathbf{L}\left(t\right)^\top \right)^\top \mathbf{X}^n, \quad \mathbf{L}\left(t_n\right) = \mathbf{V}^n \mathbf{S}^{n,\top}.
\end{align*}
The angular basis is updated by determining $\widehat{\mathbf{V}}^{n+1} \in \mathbb{R}^{N_v \times 2r}$ as an orthonormal basis of $[\mathbf{L}\left(t_{n+1}\right), \mathbf{V}^n] \in \mathbb{R}^{N_v \times 2r}$, e.g. by QR-decomposition. We store $\widehat{\mathbf{N}} = \widehat{\mathbf{V}}^{n+1,\top} \mathbf{V}^n \in \mathbb{R}^{2r \times r}$.

\textbf{\textit{S}-step}: We update the coefficient matrix from $\mathbf{S}^n \in \mathbb{R}^{r \times r}$ to $\widehat{\mathbf{S}}^{n+1} \in \mathbb{R}^{2r \times 2r}$ by solving the ODE
\begin{align*}
\dot{\widehat{\mathbf{S}}}\left(t\right) = \widehat{\mathbf{X}}^{n+1, \top} \mathbf{F} \left( \widehat{\mathbf{X}}^{n+1} \widehat{\mathbf{S}}\left(t\right) \widehat{\mathbf{V}}^{n+1,\top} \right) \widehat{\mathbf{V}}^{n+1}, \quad \widehat{\mathbf{S}}\left(t_n\right) = \widehat{\mathbf{M}} \mathbf{S}^n \widehat{\mathbf{N}}^\top.
\end{align*}

\textbf{Truncation}: We compute the singular value decomposition of $\widehat{\mathbf{S}}^{n+1} = \widehat{\mathbf{P}} \mathbf{\Sigma}\widehat {\mathbf{Q}}^\top$, where $\widehat{\mathbf{P}}, \widehat {\mathbf{Q}} \in \mathbb{R}^{2r \times 2r}$ are orthogonal matrices and $\mathbf{\Sigma} \in \mathbb{R}^{2r \times 2r}$ is the diagonal matrix containing the singular values $\sigma_1,...,\sigma_{2r}$. The new rank $r_1 \leq 2r$ is determined such that 
\begin{align*}
\left(\sum_{j=r_1+1}^{2r} \sigma_j^2\right)^{1/2} \leq \vartheta,
\end{align*}
where $\vartheta$ denotes a prescribed tolerance. We set $\mathbf{S}^{n+1} \in \mathbb{R}^{r_1\times r_1}$ to contain the $r_1$ largest singular values of $\widehat{\mathbf{S}}^{n+1}$ and $\mathbf{P}^{n+1} \in \mathbb{R}^{2r\times r_1}$ and $\mathbf{Q}^{n+1} \in \mathbb{R}^{2r\times r_1}$ to contain the first $r_1$ columns of $\widehat{\mathbf{P}}$ and $\widehat {\mathbf{Q}}$, respectively. Finally, we compute $\mathbf{X}^{n+1} = \widehat{\mathbf{X}}^{n+1}\mathbf{P}^{n+1} \in \mathbb{R}^{N_x \times r_1}$ and $\mathbf{V}^{n+1} = \widehat{\mathbf{V}}^{n+1} \mathbf{Q}^{n+1} \in \mathbb{R}^{N_v \times r_1}$.

The update of $\mathbf{f}_r^n$ after one time step is given by $\mathbf{f}_r^{n+1} = \mathbf{X}^{n+1} \mathbf{S}^{n+1} \mathbf{V}^{n+1,\top}$. Note that in the following, to simplify notation, we will write $\mathbf{f}$ instead of $\mathbf{f}_r$. 

\subsection{Dynamical low-rank approximation for the discrete optimization problem}\label{sec4.2:DLRA for the optimization}

The goal of this section consists in applying DLRA to the fully discrete gradient descent method proposed in Algorithm \ref{alg:GD fully discrete}. To this end, we reformulate the forward equations \eqref{3:time discretization forward problem} as well as the adjoint equations \eqref{3:time discretization adjoint problem} using the dynamical low-rank method with augmented BUG integrator. 

The initial low-rank factors $\mathbf{X}_m^{0,\text{for}}, \mathbf{S}_m^{0,\text{for}},$ and $\mathbf{V}_m^{0,\text{for}}$ for the forward equations \eqref{3:time discretization forward problem} are obtained by a singular value decomposition of $\mathbf{u}^0_m$, where $\mathbf{u}^0_m = \left( u^0_{jkm}\right) \in \mathbb{R}^{N_x \times N_v}$, for which the number of singular values is truncated to the initial rank $r$. In each time step, the low-rank factors $\mathbf{X}_m^{n,\text{for}}, \mathbf{S}_m^{n,\text{for}},$ and $\mathbf{V}_m^{n,\text{for}}$ are evolved according to the following scheme. 

\begin{subequations}\label{4:DLRA forward KLS}
First, we solve in parallel the equations 
\begin{align}
\mathbf{K}^{n+1,\text{for}}_m =& \  \mathbf{K}^{n,\text{for}}_m  - \Delta t \mathbf{D}^x \mathbf{K}^{n,\text{for}}_m \mathbf{V}^{n, \text{for},\top}_m  \mathbf{A}^\top \mathbf{V}^{n,\text{for}}_m+ \Delta t \frac{\Delta x}{2} \mathbf{D}^{xx} \mathbf{K}^{n,\text{for}}_m \mathbf{V}^{n,\text{for},\top}_m |\mathbf{A}|^\top \mathbf{V}^{n,\text{for}}_m\\
&+ \Delta t \diag(\sigma) \mathbf{K}^{n,\text{for}}_m \mathbf{V}^{n,\text{for},\top}_m \mathbf{E} \mathbf{V}^{n,\text{for}}_m,\nonumber\\
\mathbf{L}^{n+1,\text{for}}_m =& \ \mathbf{L}^{n,\text{for}}_m - \Delta t \mathbf{A} \mathbf{L}^{n,\text{for}}_m \mathbf{X}^{n, \text{for},\top}_m \mathbf{D}^{x,\top} \mathbf{X}^{n,\text{for}}_m + \Delta t \frac{\Delta x}{2}|\mathbf{A}| \mathbf{L}^{n,\text{for}}_m \mathbf{X}^{n, \text{for},\top}_m \mathbf{D}^{xx,\top} \mathbf{X}^{n,\text{for}}_m\\
&+ \Delta t \mathbf{E} \mathbf{L}^{n,\text{for}}_m \mathbf{X}^{n,\text{for},\top}_m \diag(\sigma) \mathbf{X}^{n,\text{for}}_m,\nonumber 
\end{align}
where $\mathbf{E} = \diag([0,-1,...,-1])$. In the next step, we perform a QR-decomposition of the augmented quantities $\left[\mathbf{K}^{n+1,\text{for}}_m, \mathbf{X}^{n,\text{for}}_m \right]$ and $\left[\mathbf{L}^{n+1,\text{for}}_m, \mathbf{V}^{n,\text{for}}_m \right]$ to obtain the augmented and time updated spatial bases $\widehat{\mathbf{X}}^{n+1,\text{for}}_m$ and angular bases $\widehat{\mathbf{V}}^{n+1,\text{for}}_m$, respectively. For the $S$-step we introduce the notation $\widetilde{\mathbf{S}}^{n,\text{for}}_m = \widehat{\mathbf{X}}^{n+1,\text{for},\top}_m \mathbf{X}^{n,\text{for}}_m \mathbf{S}^{n,\text{for}}_m \mathbf{V}^{n, \text{for},\top}_m \widehat{\mathbf{V}}^{n+1,\text{for}}_m$ and compute
\begin{align}
\widehat{\mathbf{S}}^{n+1,\text{for}}_m =& \ \widetilde{\mathbf{S}}^{n,\text{for}}_m - \Delta t \widehat{\mathbf{X}}^{n+1, \text{for},\top}_m \mathbf{D}^x \widehat{\mathbf{X}}^{n+1,\text{for}}_m \widetilde{\mathbf{S}}^{n,\text{for}}_m \widehat{\mathbf{V}}^{n+1,\text{for},\top}_m \mathbf{A}^\top \widehat{\mathbf{V}}^{n+1,\text{for}}_m\nonumber\\
&+ \Delta t \frac{\Delta x}{2} \widehat{\mathbf{X}}^{n+1,\text{for},\top}_m \mathbf{D}^{xx} \widehat{\mathbf{X}}^{n+1,\text{for}}_m \widetilde{\mathbf{S}}^{n,\text{for}}_m \widehat{\mathbf{V}}^{n+1,\text{for},\top}_m |\mathbf{A}|^\top \widehat{\mathbf{V}}^{n+1,\text{for}}_m\\
&+ \Delta t \widehat{\mathbf{X}}^{n+1,\text{for},\top}_m \diag(\sigma) \widehat{\mathbf{X}}^{n+1,\text{for}}_m \widetilde{\mathbf{S}}^{n,\text{for}}_m \widehat{\mathbf{V}}^{n+1,\text{for},\top}_m \mathbf{E} \widehat{\mathbf{V}}^{n+1,\text{for}}_m.\nonumber
\end{align}
\end{subequations}
Finally, we truncate the time-updated augmented low-rank factors for each $m=1,...,N_{IC}$ to a new rank $r_1 \leq 2r$. Then the time-updated numerical solutions of the forward problem are given by $\mathbf{u}^{n+1}_m = \mathbf{X}^{n+1,\text{for}}_m \mathbf{S}^{n+1,\text{for}}_m \mathbf{V}^{n+1,\text{for},\top}_m \in \mathbb{R}^{N_x \times N_v}$.

For the adjoint equations \eqref{3:time discretization adjoint problem} we perform a singular value decomposition of the end time solutions $\mathbf{w}_m^{N_t} = \left( w_{jkm}^{N_t}\right) \in \mathbb{R}^{N_x \times N_v}$, truncate to the prescribed initial rank $r$ and obtain the low-rank factors $\mathbf{X}_m^{N_t,\text{adj}}, \mathbf{S}_m^{N_t,\text{adj}},$ and $\mathbf{V}_m^{N_t,\text{adj}}$. Then, in each step, the low-rank factors $\mathbf{X}_m^{n,\text{adj}}, \mathbf{S}_m^{n,\text{adj}},$ and $\mathbf{V}_m^{n,\text{adj}}$ are evolved backwards in time as follows.

First, we solve in parallel the equations
\begin{align*}
\mathbf{K}^{n-1,\text{adj}}_m =& \ \mathbf{K}^{n,\text{adj}}_m + \Delta t \mathbf{D}^x \mathbf{K}^{n,\text{adj}}_m\mathbf{V}^{n,\text{adj},\top}_m \mathbf{A}^\top \mathbf{V}^{n,\text{adj}}_m+ \Delta t \frac{\Delta x}{2} \mathbf{D}^{xx} \mathbf{K}^{n,\text{adj}}_m \mathbf{V}^{n,\text{adj},\top}_m |\mathbf{A}|^\top \mathbf{V}^{n,\text{adj}}_m\\
&+ \Delta t \diag(\sigma) \mathbf{K}^{n,\text{adj}}_m \mathbf{V}^{n,\text{adj},\top}_m \mathbf{E} \mathbf{V}^{n,\text{adj}}_m,\\
\mathbf{L}^{n-1,\text{adj}}_m =& \ \mathbf{L}^{n,\text{adj}}_m + \Delta t \mathbf{A} \mathbf{L}^{n,\text{adj}}_m \mathbf{X}^{n,\text{adj},\top}_m \mathbf{D}^{x,\top} \mathbf{X}^{n,\text{adj}}_m + \Delta t \frac{\Delta x}{2} |\mathbf{A}| \mathbf{L}^{n,\text{adj}}_m \mathbf{X}^{n,\text{adj},\top}_m \mathbf{D}^{xx,\top} \mathbf{X}^{n,\text{adj}}_m\\
&+ \Delta t \mathbf{E} \mathbf{L}^{n,\text{adj}}_m \mathbf{X}^{n,\text{adj},\top}_m \diag(\sigma) \mathbf{X}^{n,\text{adj}}_m.
\end{align*}
In the next step, we perform a QR-decomposition of the augmented quantities $\left[\mathbf{K}^{n-1,\text{adj}}_m, \mathbf{X}^{n,\text{adj}}_m \right]$ and $\left[\mathbf{L}^{n-1,\text{adj}}_m, \mathbf{V}^{n,\text{adj}}_m \right]$ to obtain the augmented and time updated spatial bases $\widehat{\mathbf{X}}^{n-1,\text{adj}}_m$ and angular bases $\widehat{\mathbf{V}}^{n-1,\text{adj}}_m$, respectively. For the $S$-step we set $\widetilde{\mathbf{S}}^{n,\text{adj}}_m = \widehat{\mathbf{X}}^{n-1,\text{adj},\top}_m \mathbf{X}^{n,\text{adj}}_m \mathbf{S}^{n,\text{adj}}_m \mathbf{V}^{n,\text{adj},\top}_m \widehat{\mathbf{V}}^{n-1,\text{adj}}_m$ and compute
\begin{align*}
\widehat{\mathbf{S}}^{n-1,\text{adj}}_m =& \ \widetilde{\mathbf{S}}^{n,\text{adj}}_m + \Delta t \widehat{\mathbf{X}}^{n-1,\text{adj},\top}_m \mathbf{D}^x \widehat{\mathbf{X}}^{n-1,\text{adj}}_m \widetilde{\mathbf{S}}^{n,\text{adj}}_m \widehat{\mathbf{V}}^{n-1,\text{adj},\top}_m \mathbf{A}^\top \widehat{\mathbf{V}}^{n-1,\text{adj}}_m\\
&+ \Delta t \frac{\Delta x}{2} \widehat{\mathbf{X}}^{n-1,\text{adj},\top}_m \mathbf{D}^{xx} \widehat{\mathbf{X}}^{n-1,\text{adj}}_m \widetilde{\mathbf{S}}^{n,\text{adj}}_m \widehat{\mathbf{V}}^{n-1,\text{adj},\top}_m |\mathbf{A}|^\top \widehat{\mathbf{V}}^{n-1,\text{adj}}_m\\
&+ \Delta t \widehat{\mathbf{X}}^{n-1,\text{adj},\top}_m \diag(\sigma) \widehat{\mathbf{X}}^{n-1,\text{adj}}_m \widetilde{\mathbf{S}}^{n,\text{adj}}_m \widehat{\mathbf{V}}^{n-1,\text{adj},\top}_m \mathbf{E} \widehat{\mathbf{V}}^{n-1,\text{adj}}_m.
\end{align*}
Finally, we truncate the time-updated augmented low-rank factors for each $m=1,...,N_{IC}$ to a new rank $r_1 \leq 2r$. Then the time-updated numerical solutions of the adjoint problem are given by $\mathbf{w}^{n-1}_m = \mathbf{X}^{n-1,\text{adj}}_m \mathbf{S}^{n-1,\text{adj}}_m \mathbf{V}^{n-1,\top,\text{adj}}_m \in \mathbb{R}^{N_x \times N_v}$.

Having determined the low-rank solutions of the forward and the adjoint problems, we can use them to compute the gradient as given in \eqref{3:time discretization gradient}. For the update of the coefficients according to \eqref{2:gradient descent step} we determine the step size adaptively by a backtracking line search approach with Armijo condition similar to \cite{ScaloneEinkemmerKuschMcClarren2024} and as described in Algorithm \ref{alg:DLRA line search}. For a given step size $\eta^n$ the coefficients and the scattering coefficient are updated to $\mathbf{c}^{n+1}$ and $\bm\sigma^{n+1}$, respectively. Then the truncation error tolerance $\vartheta$ is adjusted using the given step size $\eta^n$ and the maximal absolute value of $\nabla_{\mathbf{c}} J \left(\mathbf{c}^n \right)$. We add some safety parameters $h_2$ and $h_3$ as well as a lower bound $h_1$ for the truncation tolerance. In the next step we compute the value of the goal function $J$ with the low-rank factors of the forward problem at hand and denote it with $J^n$. Then we solve the forward problem \eqref{4:DLRA forward KLS} with $\bm \sigma^{n+1}$ and the updated value of $\vartheta$ and use the obtained low-rank factors for another evaluation of the cost function $J$, for which the result is denoted by $\widebar{J}^n$. While the difference between $J^n$ and $\widebar{J}^n$ is larger than a prescribed tolerance depending on the Euclidean norm $\left\Vert \cdot \right\Vert_2$ of the gradient of the cost function, the gradient descent step size is reduced by the step size reduction factor $p$ and the procedure is repeated. 

\begin{algorithm}[h!]\small
\caption{\small Backtracking line search method for the adaptive refinement of the gradient descent step size and the DLRA rank tolerance}\label{alg:DLRA line search}
\begin{algorithmic}
\Require \hspace{0.18cm} goal function $J$,\\
\hspace{1.04cm} coefficients $\mathbf{c}^n$,\\
\hspace{1.04cm} gradient $\nabla_{\mathbf{c}} J \left( \mathbf{c}^n \right)$ computed using \eqref{3:time discretization gradient},\\
\hspace{1.04cm} low-rank factors $ \mathbf{X}^{n,\text{for}}_m, \mathbf{S}^{n,\text{for}}_m, \mathbf{V}^{n,\text{for}}_m$ of the forward problem \eqref{4:DLRA forward KLS} for $m=1,...,N_{\text{IC}}$,\\
\hspace{1.04cm} step size $\eta^n > 0$,\\
\hspace{1.04cm} rank error tolerance $\vartheta$,\\
\hspace{1.04cm} step size reduction factor $p \in \left( 0,1 \right)$,\\ 
\hspace{1.04cm} constants $h_1,h_2,h_3,h_4$ 

\Ensure refined step size $\eta^{n+1}$, refined rank error tolerance $\vartheta$, updated coefficients $\mathbf{c}^{n+1}$\\

\State Update the coefficients according to \eqref{2:gradient descent step}: $\mathbf{c}^{n+1} = \mathbf{c}^n - \eta^n \nabla_{\mathbf{c}} J \left( \mathbf{c}^n \right)$;
\State Compute $\bm{\sigma}^{n+1}$ from the updated coefficients $\mathbf{c}^{n+1}$ according to \eqref{2:sigma spline approximation discretized}; 
\State Update $\vartheta = \max \left( h_1, \min \left(h_2, \eta^n h_3 \left \Vert\nabla_{\mathbf{c}} J \left(\mathbf{c}^n \right) \right\Vert_{\infty} \right)\right)$;\\

\State Compute $J^n = J \left( \mathbf{X}^{n,\text{for}}_1 \mathbf{S}^{n,\text{for}}_1 \mathbf{V}^{n,\text{for}}_1, ..., \mathbf{X}^{n,\text{for}}_{N_{\text{IC}}} \mathbf{S}^{n,\text{for}}_{N_{\text{IC}}} \mathbf{V}^{n,\text{for}}_{N_{\text{IC}}}\right)$;
\State Compute $\widebar{\mathbf{X}}^{n,\text{for}}_m \widebar{\mathbf{S}}^{n,\text{for}}_m \widebar{\mathbf{V}}^{n,\text{for}}_m$  from \eqref{4:DLRA forward KLS} for $m=1,...,N_{\text{IC}}$ with $\bm{\sigma}^{n+1}$ and the updated $\vartheta$; 
\State Compute $\widebar{J}^n = J \left( \widebar{\mathbf{X}}^{n,\text{for}}_1 \widebar{\mathbf{S}}^{n,\text{for}}_1 \widebar{\mathbf{V}}^{n,\text{for}}_1, ..., \widebar{\mathbf{X}}^{n,\text{for}}_{N_{\text{IC}}} \widebar{\mathbf{S}}^{n,\text{for}}_{N_{\text{IC}}} \widebar{\mathbf{V}}^{n,\text{for}}_{N_{\text{IC}}}\right)$;\\

\While{$\widebar{J}^n > J^n - \eta^n h_4 \left\Vert \nabla_{\mathbf{c}} J \left(\mathbf{c}^n \right) \right\Vert^2_2$}
\State Update $\eta^{n+1} = p \eta^n$;\\

\State Update $\mathbf{c}^{n+1} = \mathbf{c}^{n+1} - \eta^{n+1} \nabla_{\mathbf{c}} J \left( \mathbf{c}^n \right)$;
\State Compute $\bm{\sigma}^{n+1}$ from the updated coefficients $\mathbf{c}^{n+1}$ according to \eqref{2:sigma spline approximation discretized}; 

\State Compute $\widebar{\mathbf{X}}^{n,\text{for}}_m \widebar{\mathbf{S}}^{n,\text{for}}_m \widebar{\mathbf{V}}^{n,\text{for}}_m$ from \eqref{4:DLRA forward KLS} for $m=1,...,N_{\text{IC}}$ with $\bm{\sigma}^{n+1}$; 
\State Compute $\widebar{J}^n = J \left( \widebar{\mathbf{X}}^{n,\text{for}}_1 \widebar{\mathbf{S}}^{n,\text{for}}_1 \widebar{\mathbf{V}}^{n,\text{for}}_1, ..., \widebar{\mathbf{X}}^{n,\text{for}}_{N_{\text{IC}}} \widebar{\mathbf{S}}^{n,\text{for}}_{N_{\text{IC}}} \widebar{\mathbf{V}}^{n,\text{for}}_{N_{\text{IC}}}\right)$;\\

\State Set $\eta^n= \eta^{n+1}$; 
\EndWhile\\

\State Update $\vartheta = \max \left( h_1, \min \left(h_2, \eta^{n+1} h_3 \left \Vert\nabla_{\mathbf{c}} J \left( \mathbf{c}^n \right) \right\Vert_{\infty} \right)\right)$;\\

\end{algorithmic}
\end{algorithm}\normalsize

\section{Numerical results}\label{sec5:Numerical results}

We consider the following test examples in both 1+1D, i.e.~one dimension of space and one dimension of angle, and 2+2D to show the computational accuracy and efficiency of the proposed low-rank scheme compared to the solution computed with the full forward and adjoint solver. 

\subsection{1+1D Cosine initial value}\label{sec5.1:1D Cosine}

For the first numerical experiment the spatial as well as the angular domain are set to $\Omega_x = \Omega_v = \left[-1,1\right]$. We consider $N_{\text{IC}} =3$ initial distributions of Cosine type of the form 
\begin{align*}
u_m \left(t=0,x \right) = 2 + \cos\left( \left( x- \frac{2m}{3}\right) \pi \right) \quad \text{ for } \quad m=1,2,3.
\end{align*}
The true and the initial B-spline coefficients for the approximation of the scattering coefficient $\sigma$ are chosen as 
\begin{align*}
\mathbf{c}_{\text{true}} = \left( 2.1,2.0,2.2\right)^\top \quad \text{ and } \quad \mathbf{c}_{\text{init}} = \left( 1.0,1.5,3.0\right)^\top. 
\end{align*}
The order of the spline basis functions is set to $3$, i.e. cubic periodic B-splines are considered. As computational parameters we use $N_x = 100$ cells in the spatial domain and $N_v =50$ moments for the approximation in the angular variable. The final time is set to $T=1.0$ and the time step size of the algorithm is chosen such that $\Delta t = \text{CFL} \cdot \Delta x$ with a $\text{CFL}$ number of $\text{CFL}=0.99$. For the low-rank computations we start with an initial rank of $r=5$ in the forward as well as in the adjoint problem. The maximal allowed value of the rank in each step shall be restricted to $20$. We begin the gradient descent method with a step size of $\eta^0 = 5 \cdot 10^5$ and a truncation error tolerance of $\vartheta = 10^{-2} \left\Vert \mathbf{\Sigma}\right\Vert_2$. For the rescaling of the gradient descent step size and the DLRA rank tolerance we use the step size reduction factor $p=0.5$ as well as the constants $h_1 = 10^{-3} \left\Vert \mathbf{\Sigma}\right\Vert_2$ for a lower bound of the rank tolerance and $h_2=0.1, h_3 = 0.1$ as safety parameters. Also $h_4 = 0.5$ is added as a safety parameter to ensure a reasonable difference between $\widebar{J}^n$ and $J^n$ in Algorithm \ref{alg:DLRA line search}. The whole gradient descent procedure is conducted until the prescribed error tolerance \texttt{errtol} $=10^{-4}$ is reached or a maximal number of iterations \texttt{maxiter} $=500$ is reached. 

\begin{figure}[t!]
    \centering
    \includegraphics[width = 0.73\linewidth]{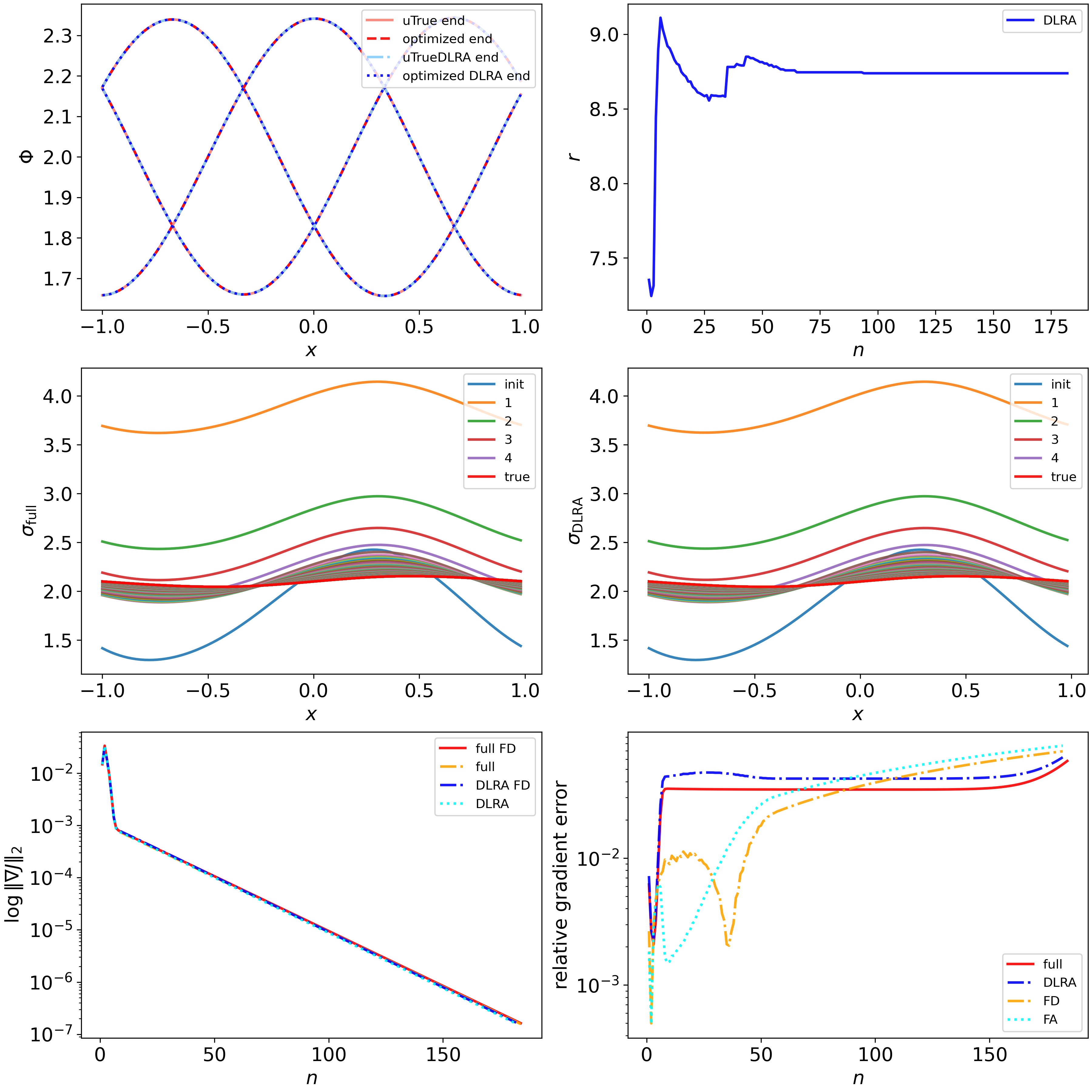}
   \caption{Top left: Numerical results for the scalar flux $\Phi$ of the 1D Cosine problem computed with the full solvers and the DLRA solvers, both with the true coefficients and with the optimization gradient descent scheme. Top right: Evolution of the averaged rank $r$ for the DLRA method. Middle row: Iterations for the reconstruction of the scattering coefficient $\sigma$ computed with both the full solvers (left) and the DLRA solvers (right). Bottom left: Evolution of the gradient norm, plotted on logarithmic scale, computed with the forward-adjoint method as well as with a forward finite difference solver for both the full and the DLRA scheme. Bottom right: Evolution of different relative gradient errors, plotted on a logarithmic scale.}
    \label{fig:Cosine}
\end{figure}

In Figure \ref{fig:Cosine} we compare the solutions of the parameter identification problem computed with the full solvers and the DLRA solvers for both the forward and the adjoint equations. We plot three curves corresponding to the different initial conditions of the scalar flux $\Phi = \frac{1}{\sqrt{2}}\left\langle f \right \rangle_v$ and display the results that are obtained from the true coefficients $\left(\text{uTrue end}\right)$ and at the end of the gradient descent optimization procedure $\left(\text{optimized end}\right)$, evaluated with both the full and the DLRA solver. We observe that the DLRA solution captures well the behavior of the full solution and that they both approach the solutions computed with the true coefficients. In addition, the parameter reconstruction inverse problem for determining $\sigma$ is resolved accurately with both solvers. It can be seen in the middle row that beginning with $\bm{\sigma}_{\text{init}}$ both the full and the DLRA method converge to the true solution $\bm{\sigma}_{\text{true}}$. Note that in the legend of both plots only the initial configuration (init) as well as the first 4 steps of the iteration converging towards the true scattering coefficient (true) are shown. Further, on the top right the evolution of the rank $r$ is depicted, where we have averaged the ranks of the forward equations computed with the different initial conditions to obtain $r_{\text{for}}$ and the ranks of the adjoint equations computed with the different initial conditions to obtain $r_{\text{adj}}$ and finally set $r = \frac{1}{2} \left(r_{\text{for}} + r_{\text{adj}} \right)$. We observe that in the beginning the averaged rank decreases as the initial rank was chosen larger than required. From then on, we observe a relatively monotonous increase until it stays at approximately $r=9$. This evolution of the rank reflects the fact that in the beginning of the optimization the error tolerance $\vartheta$ is chosen quite large as the computed solution is still comparably far away from the true solution. As the optimization algorithm approaches the true coefficients, the DLRA rank tolerance $\vartheta$ is decreased, resulting in a higher averaged rank. In the bottom row the evolution of $\left\Vert \nabla J \right\Vert_2$ is plotted on a logarithmic scale. For both the full and the DLRA scheme, it is computed with the presented forward-adjoint method as well as with a forward finite difference scheme (FD) for comparison. It can be observed that the results of each computation lead to similar results. In addition, it can be observed that $\left\Vert \nabla J \right\Vert_2$ mainly decreases with the number of performed iterations, indicating that the searched-for solution is steadily approached. On the bottom left different relative gradient errors are displayed. We compare the results obtained from the forward-adjoint method with the results obtained from the forward difference scheme for both the full solvers (full) and the DLRA solvers (DLRA). Further, the difference between the forward difference scheme with full and DLRA solvers (FD) and the values obtained from the forward-adjoint method with full and DLRA solvers (FA) are depicted. We observe that the low-rank approximation adds an error that is similar to the difference between computing the gradient with forward differences and computing it using the adjoint.

For the considered setup, the computational benefit of the DLRA method compared to the solution of the full problem is significant. Written in Julia v1.11 and run on a MacBook Pro with M1 chip, the run time decreases by a factor of approximately $1.52$ from $34.49$ seconds to $22.66$ seconds while retaining the accuracy of the computed results. Concerning the memory costs, the solutions of the forward problem and of the adjoint problem have to be stored in order to compute the gradient. For each initial condition, the storage of the solution of the forward problem corresponds to a memory cost of $8 \left(N_t+1 \right) N_x N_v$, which for the DLRA method can be lowered to $8\left(N_t+1 \right) \left(rN_x+rN_v+r^2 \right)$, where $r$ is the maximal averaged rank in the simulation.

\subsection{1+1D Gaussian initial value}

In a second test example, we set $\Omega_x = \left[0,10\right]$ for the spatial and $\Omega_v = \left[-1,1\right]$ for the angular domain. We consider $N_{\text{IC}} =5$ Gaussian initial distributions of the form 
\begin{align*}
u_m \left(t=0,x \right) =\max \left(10^{-8}, \frac{1}{\sqrt{2\pi \sigma_{\text{IC}}^2}} \exp{\left(-\frac{\left(x-x_0\right)^2}{2 \sigma_{\text{IC}}^2}\right)} \right)  \quad \text{ for } \quad m=1,2,3,4,5,
\end{align*}
which are centered around equidistantly distributed $x_0$ and extended periodically to the domain $\Omega_x$. The standard deviation is set to the constant value $\sigma_{\text{IC}}=0.8$. The true and the initial spline coefficients for the approximation of the scattering coefficient $\sigma$ are chosen as 
\begin{align*}
\mathbf{c}_{\text{true}} = \left( 2.1,2.0,2.2,2.0,1.9\right)^\top \quad \text{ and } \quad \mathbf{c}_{\text{init}} = \left( 2.8,1.5,3.0,2.1,1.2\right)^\top. 
\end{align*}
All other settings and computational parameters remain unchanged from the previous test example. 

\begin{figure}[t!]
    \centering
    \includegraphics[width = 0.73\linewidth]{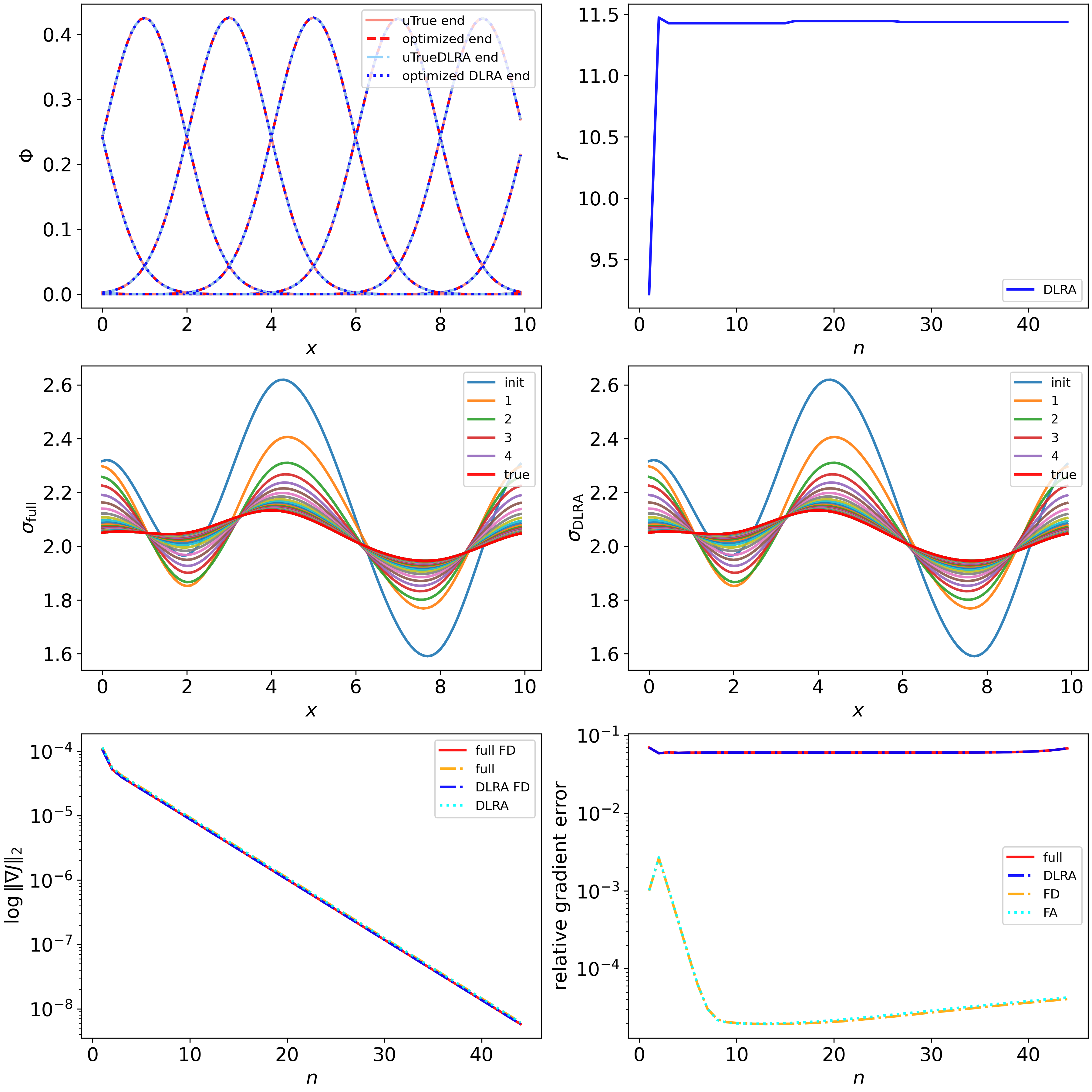}
   \caption{Top left: Numerical results for the scalar flux $\Phi$ of the 1D Gauss problem computed with the full solvers and the DLRA solvers, both with the true coefficients and with the optimization gradient descent scheme. Top right: Evolution of the averaged rank $r$ for the DLRA method. Middle row: Iterations for the reconstruction of the scattering coefficient $\sigma$ computed with both the full solvers (left) and the DLRA solvers (right). Bottom left: Evolution of the gradient norm, plotted on logarithmic scale, computed with the forward-adjoint method as well as with a forward difference solver for both the full and the DLRA scheme. Bottom right: Evolution of different relative gradient errors, plotted on logarithmic scale.}
    \label{fig:Gauss}
\end{figure}

In Figure \ref{fig:Gauss} we compare the solutions of the parameter identification problem computed with the full solvers and the DLRA solvers for both the forward and the adjoint equations. We plot five curves corresponding to the different initial conditions of the scalar flux $\Phi = \frac{1}{\sqrt{2}}\left\langle f \right \rangle_v$, computed with the true coefficients $\left(\text{u True end}\right)$ and at the end of the optimization procedure $\left(\text{optimized end}\right)$, evaluated with both the full and the DLRA solver. Again we observe that the DLRA solution captures well the behavior of the full solution and that they both approach the solutions computed with the true coefficients. For the reconstruction of the scattering coefficient $\sigma$ it can be seen that beginning with $\bm{\sigma}_{\text{init}}$ both the full and the DLRA method converge to the true solution $\bm{\sigma}_{\text{true}}$. The averaged rank $r$ first decreases as the initial rank was chosen larger than required. From then on, we observe the expected relatively monotonous increase until it stagnates at a value of approximately $r=11.5$. On the bottom left a monotonous decrease of $\left\Vert \nabla J \right\Vert_2$, computed with the forward-adjoint method, can be observed. We can also see (bottom right) that the error of the DLRA is much smaller than the difference between forward difference and the adjoint method.

Written in Julia v1.11 and run on a MacBook Pro with M1 chip, the computational time of the DLRA method compared to the solution of the full problem decreases by a factor of approximately $1.13$ from $4.09$ seconds to $3.61$ seconds.

\subsection{2+2D Cosine initial value \label{sec:4d-example}}

The computational benefits of the presented DLRA scheme are significantly enhanced in higher dimensions. To demonstrate this, we extend the numerical considerations to a 2D spatial and a 2D angular setting. The corresponding forward equations for $m=1,...,N_{\text{IC}}$ become 
\begin{align*}
\begin{cases}
\partial_t f_m \left(t,\mathbf{x}, \mathbf{\Omega}\right) + \mathbf{\Omega} \cdot \nabla_{\mathbf{x}} f_m \left(t,\mathbf{x}, \mathbf{\Omega} \right) &= \sigma \left(\mathbf{x} \right) \left(\frac{1}{|\Omega_{\mathbf{\Omega}}|}\langle f_m \left(t,\mathbf{x}, \mathbf{\Omega}\right)) \rangle_{\mathbf{\Omega}} - f_m \left(t,\mathbf{x}, \mathbf{\Omega}\right) \right),\\
f_m\left(t=0,\mathbf{x}, \mathbf{\Omega}\right) &= f_{\text{in},m}\left(\mathbf{x}, \mathbf{\Omega}\right).
\end{cases}
\end{align*}
For our numerical experiments we set $\mathbf{x} = \left( x,y \right) \in \left[-1,1 \right] \times \left[-1,1 \right]$ and $\mathbf{\Omega} \in \mathcal{S}^2 $. For simplicity, we consider $\sigma \left( \mathbf{x} \right) = \sigma \left(x \right)$, i.e. we consider a material that is homogeneous in the y-direction.  As initial conditions we prescribe a Cosine distributions of the form
\begin{align*}
f_m \left(t=0,\mathbf{x}, \mathbf{\Omega} \right) = 2 + \cos\left( \left( x- \frac{2m}{3}\right) \pi \right) \sin\left( \left( y- \frac{2m}{3}\right) \pi \right) \quad \text{ for } \quad m=1,2,3.
\end{align*}
The true and the initial B-spline coefficients for the approximation of the scattering coefficient $\sigma (x)$ are chosen equivalently to the 1+1D Cosine test example presented in Section \ref{sec5.1:1D Cosine}. We perform our computations on a spatial grid with $50$ points in $x$ and $50$ points in $y$. For
the angular basis we use spherical harmonics and and a polynomial degree of $n_{\text{PN}} = 19$, corresponding to $400$ expansion coefficients in angle. The end time is set to $T=1.0$ and the time step size of the algorithm is chosen such that $\Delta t = \text{CFL} \cdot \Delta x$ with a $\text{CFL}$ number of $\text{CFL}=0.7$. For the low-rank computations we start with an initial rank of $r=20$ in the forward as well as in the adjoint problem. The maximal allowed value of the rank in each step is restricted to $100$. The initial gradient descent step size is chosen as $\eta^0 = 50$, while an initial truncation error tolerance of $\vartheta = 10^{-5} \left\Vert \mathbf{\Sigma}\right\Vert_2$ is used. For the rescaling of the gradient descent step size and the DLRA rank tolerance we use the step size reduction factor $p=0.5$ as well as the constants $h_1 = 10^{-5} \left\Vert \mathbf{\Sigma}\right\Vert_2$ for a lower bound of the rank tolerance and $h_2=0.1, h_3 = 0.5, h_4 = 0.5$ as safety parameters. The gradient descent procedure terminates when a prescribed error tolerance of \texttt{errtol} $=10^{-6}$ is reached or the maximal number of iterations \texttt{maxiter} $=500$ is reached.

\begin{figure}[t!]
    \centering
    \includegraphics[width = 0.55\linewidth]{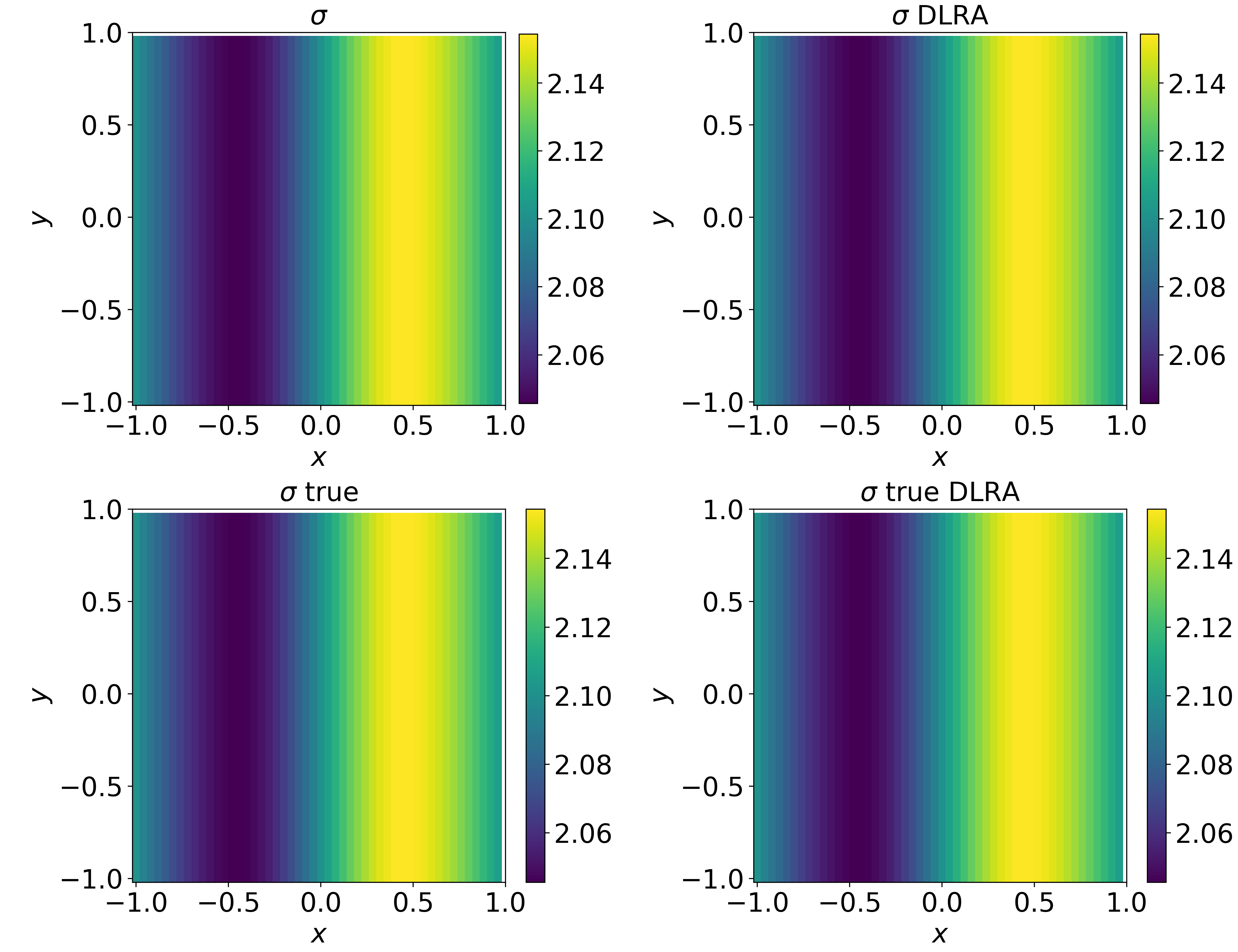}
   \caption{Numerical results for the reconstruction of the scattering coefficient $\sigma (x)$, computed with the full solvers (top left) and the DLRA solvers (top right). The bottom row shows the results obtained with the true scattering coefficient, computed with both schemes.}
    \label{fig:Cosine2D}
\end{figure}

\begin{figure}[t!]
    \centering
    \includegraphics[width = 0.7\linewidth]{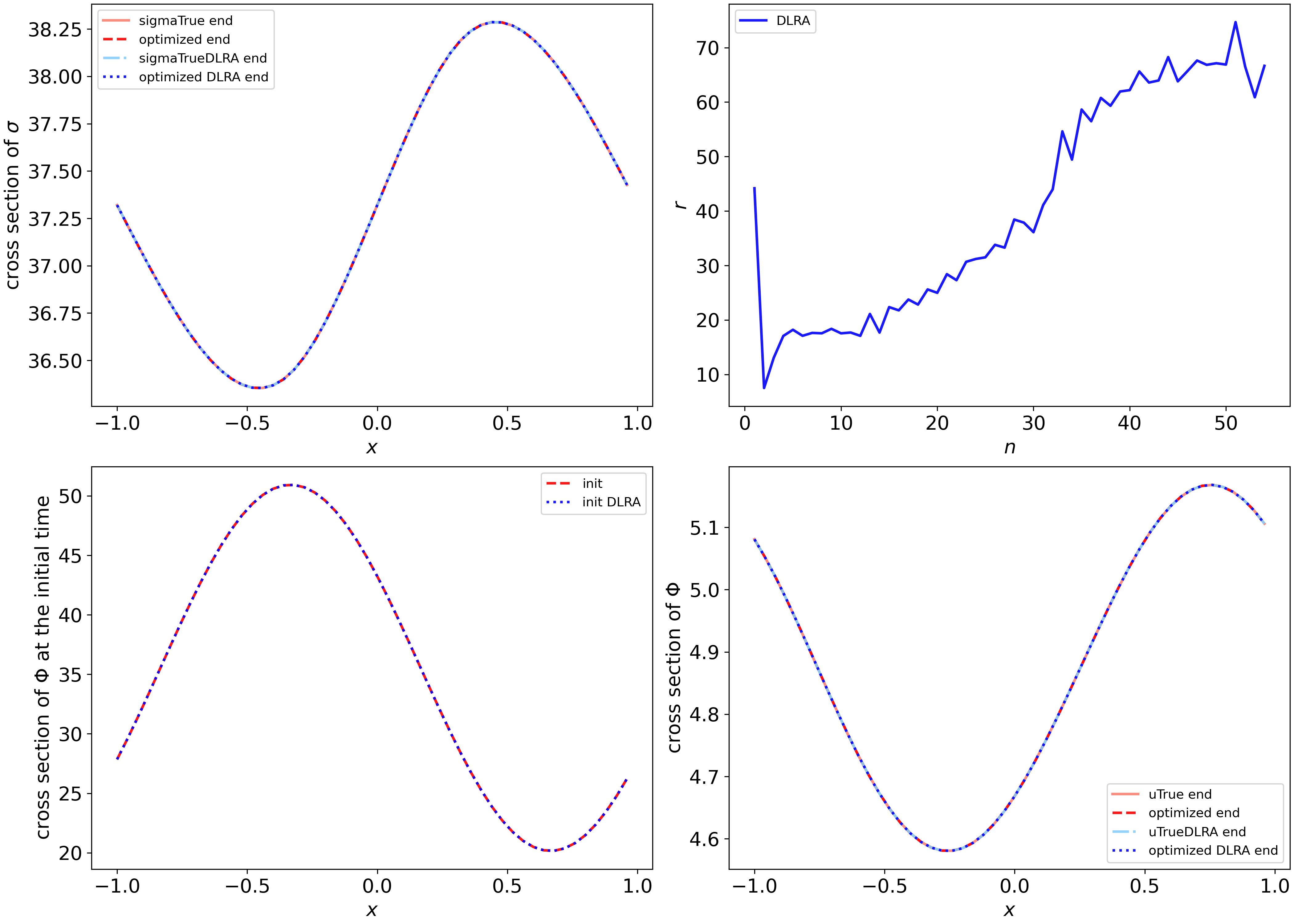}
    \caption{Top left: Cross section through the middle of the spatial domain in $y$ of the results for the reconstruction of the scattering coefficient $\sigma (x)$ obtained with the full solvers and the DLRA solvers at the end of the optimization and computed with the corresponding ground truth parameter with both solvers. Top right: Evolution of the averaged rank $r$ for the DLRA method. Bottom row: Numerical results for the cross sections through the middle of the spatial domain in $y$ of the scalar flux $\Phi$ at the initial configuration (left) and at the end of the optimization (right).}
    \label{fig:Cosine2D-allplots}
\end{figure}

In Figure \ref{fig:Cosine2D} the results from the parameter reconstruction problem for determining the scattering coefficient $\sigma = \sigma (x)$ are displayed. We observe that the DLRA scheme accurately reconstructs the scattering coefficient compared with the results obtained from the full scheme and compared with the results obtained when, for comparison, we insert the true $\sigma$ into both schemes. The top row of Figure \ref{fig:Cosine2D-allplots} shows a cut through the results obtained for the scattering coefficient as well as the evolution of the averaged rank in time. We observe that the rank increases as the optimization procedure approaches the true solution, as expected. In the bottom row cross sections of the scalar flux $\Phi = \int_{\mathcal{S}^2} f \left( t, \mathbf{x}, \mathbf{\Omega} \right) \mathrm{d}\mathbf{\Omega}$ are presented. The left plot shows a cross section of the initial distribution, computed with the full and the DLRA solvers. In the right plot the solutions obtained with the full solvers and the DLRA solvers at the end of the optimization procedure are compared with the results obtained with the true coefficient, computed with both schemes. One observes that also the result for the scalar flux is accurately captured with the DLRA scheme. 

Written in Julia v1.11 and run on a MacBook Pro with M1 chip, the computational time of the DLRA method compared to the solution of the full problem decreases by a factor of approximately $4.4$ from $6425$ seconds to $1466$ seconds. The memory required to store a single snapshot of the full numerical solution is $8$ MB (for the parameters chosen here) vs. between $200$ kB ($r=10$) and $1.4$ MB ($r=70$) for the DLRA. Even in the worst case this is a reduction by more than a factor of $5$. These results underline the ability of the DLRA scheme to obtain accurate and efficiency solutions of the parameter identification problem.

\section{Conclusion and outlook}\label{sec6:Conclusion and outlook}

We have presented a fully discrete DLRA scheme for the reconstruction of the scattering parameter in the radiative transfer equation making use of a PDE constrained optimization procedure. For a further enhancement of its computational advantages compared to a standard full solver, the step size of the gradient descent method is determined adaptively in each step and the allowed DLRA rank tolerance is adjusted accordingly. As we have seen, this leads to an efficient and accurate numerical DLRA scheme. 

In contrast to the ``optimize then low-rank" strategy we have pursued, one could consider a ``low-rank then optimize" approach. However, in this case it is not so clear if a numerically efficient form of the adjoint equations can be obtained, as the low-rank equations are highly nonlinear. It is thus not clear what, if any, advantage such an approach would provide. Another direction of potential future research could be to extend the stability analysis done for the forward problem (see, e.g., \cite{BaumannEinkemmerKlingenbergKusch2024}) to the entire parameter identification problem. This would be particularly interesting in the case where some noise is added to the input data.

\section*{Acknowledgements}

The authors thank Kui Ren for helpful and inspiring discussions. Lena Baumann acknowledges support by the W\"urzburg Mathematics Center for Communication and Interaction (WMCCI) as well as the Stiftung der Deutschen Wirtschaft (Foundation of German Business). 

\newpage
\bibliographystyle{abbrv}
\bibliography{references}

\end{document}